\documentclass[11pt,a4paper]{article}   

\usepackage[centertags,leqno]{amsmath}
\usepackage{amssymb,amsfonts,theorem,stmaryrd}
\usepackage[all]{xy}

\usepackage{amsmath}
\usepackage{amssymb,amsfonts,theorem,stmaryrd}
\usepackage[all]{xy}

\allowdisplaybreaks

\numberwithin{equation}{section}

\theorembodyfont{\slshape}
\newtheorem{thm}{Theorem}[section]
\newtheorem{prop}[thm]{Proposition}

\newtheorem{lemma}[thm]{Lemma}
\newtheorem{cor}[thm]{Corollary}
\theorembodyfont{\rmfamily} 

\newtheorem{rem}[thm]{Remark}

\newenvironment{pf}{\textit{Proof.}}{\hfill$\boxbox$}

\newcommand{\R}{\mathbb{R}}
\newcommand{\N}{\mathbb{N}}
\newcommand{\Z}{\mathbb{Z}}
\newcommand{\C}{\mathbb{C}}
\renewcommand{\H}{\mathbb{H}}

\newcommand{\E}{\mathcal{E}}

\newcommand{\Sph}{\mathbf{S}}

\newcommand{\Ric}{\operatorname{Ric}}

\newcommand{\dom}{\operatorname{dom}}
\newcommand{\spec}{\operatorname{spec}}

\newcommand{\im}{\operatorname{im}}

\renewcommand{\min}{min}
\renewcommand{\max}{max}

\renewcommand{\sl}{\mathfrak{sl}}
\newcommand{\su}{\mathfrak{su}}
\newcommand{\SL}{\operatorname{SL}}

\newcommand{\SU}{\operatorname{SU}}

\newcommand{\End}{\operatorname{End}}

\newcommand{\hol}{\operatorname{hol}}
\newcommand{\id}{\operatorname{id}}

\newcommand{\tr}{\operatorname{tr}}

\newcommand{\rank}{\operatorname{rank}}

\newcommand{\Defo}{\operatorname{Def}}

\newfont{\hugemath}{cmsy10 scaled 3000}

\begin{document}
\title{The deformation theory of hyperbolic cone-3-manifolds with cone-angles less than $2\pi$.}
\author{Hartmut Wei{\ss} \\ LMU M\"unchen}

\maketitle

\begin{abstract}
We develop the deformation theory of hyperbolic cone-3-manifolds with cone-angles less than $2\pi$, i.e.~contained in the interval $(0,2\pi)$. In the present paper we focus on deformations keeping the topological type of the cone-manifold fixed. We prove local rigidity for such structures. This gives a positive answer to a question of A.~Casson.
\end{abstract}

\tableofcontents

\section{Introduction}
Let $X$ be a closed, orientable hyperbolic cone-3-manifold. Recall that $X$ is a path metric space homeomorphic to a closed, orientable 3-manifold with certain local models prescribed, cf.~\cite{BLP} or \cite {CHK}. More precisely, for $x \in X$ the metric ball $B_\varepsilon(x) \subset X$ is required to be isometric to a truncated hyperbolic cone over a space $\mathbf{S}_x$, which in turn is required to be a spherical cone-surface homeomorphic to the 2-sphere, see below. The space $\mathbf{S}_x$ is called the {\em link} of $x$. Let $M$ denote the open subset of $X$ consisting of those points $x \in X$ with the property that $\mathbf S_x$ is actually isometric to $\mathbf S^2$, the standard smooth round 2-sphere. This subset $M$ carries a smooth but typically incomplete hyperbolic metric $g^{TM}$ and is called the {\em smooth part} of $X$. The complement $\Sigma=X \setminus M$ is called the {\em singular locus} of $X$. 

Now a spherical cone-surface $S$ is again a path metric space, homeomorphic to a surface in this case, with certain local models prescribed: A metric ball $B_\varepsilon(x) \subset S$ is required to be isometric to a truncated spherical cone over $\mathbf{S}^1_{\alpha_x}=\R / \alpha_x \Z$ for some number $\alpha_x >0$. The number $\alpha_x>0$ is called the {\em cone-angle} at $x$. The smooth part $N$ consisting of points $x \in S$ with $\alpha_x = 2\pi$ carries a smooth, but typically incomplete spherical metric $g^{TN}$ and the singular locus is just a finite collection of points $\{p_1, \ldots,p_m\}$, also called {\em cone-points} in the following. This concludes the description of the local structure of a hyperbolic cone-3-manifold. 

From the description of the local structure of a hyperbolic cone-3-manifold it is evident that the singular locus $\Sigma \subset X$ is a geodesic graph. Let $e_1, \ldots, e_N$ denote the edges and $v_1, \ldots, v_k$ the vertices contained in $\Sigma$. Note that $\Sigma$ may be disconnnected and that some edges $e_i$ may be closed singular geodesics. For each vertex $v_j$ let $m_j$ denote the number of egdes meeting at $v_j$, or equivalently, the number of cone-points contained in $\mathbf S_j$, the link of $v_j$. To each edge $e_i$ we attach a number $\alpha_i >0$, the {\em cone-angle} along $e_i$, in the following way: For each point $x$ in the interior of $e_i$ the link $\mathbf S_x$ is isometric to a spherical suspension $\mathbf S_i=\mathbf S^2 (\alpha_i,\alpha_i)$, i.e.~a spherical cone-manifold structure on $S^2$ with two cone-points and both cone-angles equal to $\alpha_i$. 

If the cone-angles are assumed to be less than or equal to $\pi$, then the list of possible links is short: $\mathbf S_x$ may either be the smooth round 2-sphere $\mathbf{S}^2$, a spherical suspension $\mathbf S^2 (\alpha,\alpha)$ as above or a cone-surface of type $\mathbf{S}^2(\alpha,\beta,\gamma)$, i.e.~the double of a spherical triangle with interior angles $\alpha/2,\beta/2$ and $\gamma/2$. All these cone-surfaces are rigid in the sense that their isometry type is determined by the cone-angles. On the other hand, if the cone-angles are allowed to lie in the interval $(0,2\pi)$, then there is a much larger choice of possible links. Moreover, these links are in general more flexible, i.e.~their isometry type is not determined by the cone-angles alone. More precisely, as a consequence of \cite{Tro} and \cite{LT}, see also \cite{MaW}, one has that for $m \geq 3$ the space of spherical cone-manifold structures on $(S^2,\{p_1, \ldots,p_m\})$ is locally parametrized by $\mathcal{T}_{0,m} \times (0,2\pi)^m$, where $\mathcal{T}_{0,m}$ is the Teichm\"uller space of the $m$-times punctured sphere. One part of the original motivation for this present work was to understand how this additional flexibility of the links affects the deformation theory of hyperbolic cone-3-manifolds.

Now let $C_{-1}(X,\Sigma)$ denote the space of hyperbolic cone-manifold structures on $(X,\Sigma)$, i.e.~hyperbolic cone-manifold structures on $X$ with singular locus precisely given by $\Sigma$. The pair $(X,\Sigma)$ is called the {\em topological type} of the cone-manifold structure in question. The space $C_{-1}(X,\Sigma)$ carries a topology such that the map
$$
\alpha: C_{-1}(X,\Sigma) \rightarrow \R^N_+
$$
mapping a hyperbolic cone-manifold structure to the vector of cone-angles $\alpha=(\alpha_1, \ldots,\alpha_N) \in \R^N_+$ is continuous. 

In \cite{HK}, C.D.~Hodgson and S.P.~Kerckhoff showed that if the cone-angles are less than $2\pi$ and the singular locus is assumed to contain no vertices, i.e.~$\Sigma$ is a disjoint union of circles, then $\alpha$ is a local homeomorphism at the given structure.

After the appearance of the results of Hodgson and Kerckhoff, A.~Casson asked in a conference talk, if $\alpha$ is a local homeomorphism in the case that vertices are present and the cone-angles are less than $2\pi$. At the same time he presented a counterexample in the case that some of the cone-angles are larger than $2\pi$. Moreover he asked, if $\alpha$ is at least always open.

In \cite{Wei}, the author showed that if the cone-angles are less than or equal to $\pi$ and $\Sigma$ is allowed to contain vertices, then again $\alpha$ is a local homeomorphism at the given structure. Under the same condition on the cone-angles the spherical case is also treated in \cite{Wei}, whereas the Euclidean case is treated in \cite{PW}. In both of these cases an additional non-degeneracy condition has to be imposed; the Euclidean case involves deforming into nearby hyperbolic and spherical structures. For details we refer the reader to \cite{Wei} and \cite{PW}.

The aim of this present work is to bridge the gap between the results contained in \cite{HK} and \cite{Wei}, namely we prove the following:
\begin{thm}[Local Rigidity]\label{localrigidity}
Let $X$ be a hyperbolic cone-3-manifold with cone-angles less than $2\pi$. Then the map 
$$
\alpha=(\alpha_1, \ldots, \alpha_N): C_{-1}(X,\Sigma) \rightarrow \R_+^N
$$ 
is a local homeomorphism at the given structure.
\end{thm}
In fact, we show slightly more, namely that the deformation space $\Defo(M)$ of incomplete hyperbolic structures on $M$, the smooth part of $X$, is smooth near the given strucuture and of complex dimension $N + \sum_{j=1}^k (m_j-3)$. We identify the deformations which correspond to cone-manifold structures preserving $(X,\Sigma)$, i.e.~$C_{-1}(X,\Sigma) \subset \Defo(M)$, together with a good local parametrization of these deformations. This essentially yields Theorem \ref{localrigidity}, which in particular gives a positive answer to the above mentioned question of Casson.

What remains to be done is to find a good local parametrization for the whole of $\Defo(M)$ and a geometric description of those deformations which are transverse to $C_{-1}(X,\Sigma)$. We will return to this issue in a joint work with G.~Montcouquiol, cf.~\cite{MoW}.

Let us finally mention that R.~Mazzeo and G.~Montcouquiol have developed an alternative approach to these questions using the deformation theory of the Einstein equation, cf.~\cite{MaM} and \cite{Mo}.

The author would like to thank Steve Kerckhoff and Rafe Mazzeo for useful conversations during the preparation of this article.

\section{Analysis on manifolds with conical singularities}\label{analysis}

\subsection{$L^2$-cohomology}

Let $(M^m,g^{TM})$ be a Riemannian manifold and $(\E,\nabla^\E,h^\E)$ a flat vector bundle over $M$, where we do not assume $h^\E$ to be $\nabla^\E$-parallel. In our main instance of such a situation, $M$ will be the smooth part of a hyperbolic cone-3-manifold and $\E$ the flat bundle of infinitesimal isometries equipped with its canonical metric.

Let $\Omega^\bullet(M;\E)$ denote the smooth differential forms on $M$ with
values in $\E$ and let $\Omega_{L^2}^\bullet(M;\E)=\{ \omega \in\Omega^\bullet(M;\E): \omega \in L^2, d^\E\omega \in L^2 \}$, where $d^\E$ is the exterior differential associated with the flat connection $\nabla^\E$. Clearly $\Omega_{L^2}^\bullet(M;\E)$ is a complex with differential $d^\E$ and its cohomology $H^\bullet_{L^2}(M;\E)$ is by definition the $L^2$-cohomology of $M$ with values in $\E$.

If we view $d^\E$ as an unbounded operator acting on compactly supported smooth forms, we consider the closed extensions $d^\E_{max}$ and $d^\E_{min}$. Let $H^\bullet_{max}$ denote the cohomology of the $d^\E_{max}$-complex and $H^{\bullet}_{min}$ the cohomology of the $d^\E_{min}$-complex. If the $L^2$-Stokes theorem holds for $\E$-valued forms on $M$, then $d^\E_{max}=d^\E_{min}$. Let
$$
\Delta_{d^\E_{max}} = \delta^\E_{min}d^\E_{max} + d^\E_{max}\delta^\E_{min}
$$
and
$$
\Delta_{d^\E_{min}} = \delta^\E_{max}d^\E_{min} + d^\E_{min}\delta^\E_{max}.
$$
Both $\Delta_{d^\E_{max}}$ and $\Delta_{d^\E_{min}}$ are selfadjoint extensions of $\Delta_{d^\E}=\delta^\E d^\E + d^\E\delta^\E$, the latter considered as acting on compactly supported smooth forms. We denote by 
$$
\mathcal H^\bullet_{max} = \{ \omega \in \dom \Delta^\bullet_{d^\E_{max}} : d^\E_{max}(\omega) = \delta^\E_{min}(\omega) = 0 \}
$$
the $d^\E_{max}$-harmonic (or $L^2$-harmonic) and by
$$
\mathcal H^\bullet_{min} = \{ \omega \in \dom \Delta^\bullet_{d^\E_{max}} : d^\E_{min}(\omega)= \delta^\E_{max}(\omega)=0 \}
$$
the $d^\E_{min}$-harmonic forms. Note that $(d^\E_{max})^*=\delta^\E_{min}$ and $(d^\E_{min})^*=\delta^\E_{max}$.
The $d^\E_{max}$-complex computes the $L^2$-cohomology:
\begin{thm}
The inclusion of the subcomplex $\Omega_{L^2}^\bullet(M;\E) \subset \dom (d^\E_{\max})^\bullet$ induces
an isomorphism in cohomology $H^\bullet_{L^2}(M;\E) \cong  H^\bullet_{\max}$.
\end{thm}
Furthermore, one has the following $L^2$-Hodge theorem:
\begin{thm}\label{L2Hodge}
There is an orthogonal decomposition
$$
\Gamma_{L^2}(M;\Lambda^\bullet T^*M\otimes\E)=\mathcal{H}_{\max}^\bullet\oplus\overline{\im (d^\E_{\max})^{\bullet-1}}\oplus\overline{\im (\delta^\E_{\min})^\bullet}\,.
$$
\end{thm}
\begin{cor}\label{CorL2Hodge}
If the range of $d^\E_{max}$ is closed, then there is an orthogonal decomposition
$$
\Gamma_{L^2}(M;\Lambda^\bullet T^*M\otimes\E)=\mathcal{H}_{\max}^\bullet\oplus\im (d^\E_{\max})^{\bullet-1}\oplus\im (\delta^\E_{\min})^\bullet
$$
and there is a canonical isomorphism $\mathcal{H}^{\bullet}_{max} \cong H_{max}^\bullet$.
\end{cor}
If in addition $b^\E \in \Gamma(M;\E^* \otimes \E^*)$ is a fiberwise non-degenerate symmetric bilinear form which is $\nabla^\E$-parallel, then the map $\alpha \mapsto ((h^\E)^{-1} \circ b^\E) \star \alpha$ induces an isomorphism
$\mathcal H^k_{max} \cong \mathcal H^{m-k}_{min}$
and the bilinear pairing
\begin{align*}
H^k_{max} \times H^{m-k}_{min} &\rightarrow \R\\
([\alpha],[\beta]) &\mapsto \int_M b^\E(\alpha \wedge \beta),
\end{align*}
is non-degenerate, where
$\wedge: \Lambda^kT^*M \otimes \E \times \Lambda^lT^*M \otimes \E \rightarrow  \Lambda^{k+l}T^*M \otimes \E \otimes \E$.
In our main instance, $b^\E_x$ will be the Killing form on $\E_x \cong \mathfrak{sl}_2(\C)$.
\begin{thm}
If the $L^2$-Stokes theorem holds for $\E$-valued forms on $M$, then Poincar\'e duality holds for $L^2$-cohomology, i.e.~the bilinear pairing
\begin{align*}
H^k_{L^2}(M;\E) \times H^{m-k}_{L^2}(M;\E) &\rightarrow \R\\
([\alpha],[\beta]) &\mapsto \int_M b^\E(\alpha \wedge \beta),
\end{align*}
is non-degenerate.
\end{thm}

Most of these results are due to J.~Cheeger and can be found in \cite{Ch1}, see also the references in \cite{Wei}. A more recent reference is \cite{BL}, which we found especially useful. 

\subsubsection{The Hodge-Dirac operator }
First we calculate the Hodge-Dirac operator near a vertex $v$: Let $(N^n,g^{TN})$ be the spherical link of $v$.  Then the hyperbolic metric has the form
$$
g_{hyp}=dr^2 + \sinh(r)^2 g^{TN}
$$
on $U_\varepsilon(v)$. We identify triples $[\phi_0, \phi_1, \phi_2]^t$ with $\phi_i \in \Gamma(U_\varepsilon(v);\pi_N^* \Lambda^iT^*N)$ with either even or odd forms on $U_\varepsilon(v)$ as follows:
$$
\xi^{ev} = \sinh(r)^{-1} \bigl( \phi_0 + \sinh(r) \phi_1 \wedge dr + \sinh(r)^2 \phi_2\bigr)
$$
and
$$
\xi^{odd} = \sinh(r)^{-1} \bigl(\phi_0 \wedge dr + \sinh(r) \phi_1 + \sinh(r)^2 \phi_2 \wedge dr\bigr).
$$
This induces isometric isomorphisms
$$
L^2((0,\varepsilon), \Gamma_{L^2}(N; \Lambda^\bullet T^*N)) \cong \Gamma_{L^2}(U_\varepsilon(v); \Lambda^{ev}T^*M)
$$
and
$$
L^2((0,\varepsilon), \Gamma_{L^2}(N; \Lambda^\bullet T^*N)) \cong \Gamma_{L^2}(U_\varepsilon(v); \Lambda^{odd}T^*M).
$$
Then with respect to these identifications
$$
D^{ev} \begin{bmatrix} \phi_0 \\ \phi_1 \\ \phi_2 \end{bmatrix}= \left(\partial_r + \frac{1}{\sinh(r)} \left(D_N + \cosh(r) 
\begin{bmatrix} -1 & 0 &0 \\
 0& 0 &0 \\
0 & 0& 1\end{bmatrix} \right)\right)\begin{bmatrix} \phi_0 \\ \phi_1 \\ \phi_2 \end{bmatrix}.
$$
The corresponding
operator in the Euclidean situation will serve as a model operator and
is given by
$$
P_B=\partial_r +
\frac{1}{r} B  
$$
with
$$
B =  D_N +  
\begin{bmatrix} -1 & 0 & 0 \\
 0& 0 &0 \\
 0& 0 & 1\end{bmatrix}
$$
acting on triples $[\phi_0, \phi_1, \phi_2]^t$ as above.
Similarly we get
$$
D^{odd}\begin{bmatrix} \phi_0 \\ \phi_1 \\ \phi_2 \end{bmatrix} = \left(-\partial_r + \frac{1}{\sinh(r)} \left(D_N + \cosh(r) 
\begin{bmatrix} -1 & 0 &0 \\
 0& 0 &0 \\
0 & 0& 1\end{bmatrix} \right)\right)\begin{bmatrix} \phi_0 \\ \phi_1 \\ \phi_2 \end{bmatrix}
$$
and the model operator is given by
$$
P_B^t = - \partial_r + \frac{1}{r} B = - P_{-B}
$$
acting on triples $[\phi_0, \phi_1, \phi_2]^t$ as above. Note that $B$ is a
symmetric first order differential operator on $N$ of conic type. It was shown in \cite{Wei} that any self-adjoint extension of $B$ has discrete spectrum.

For purposes of exposition we give a detailed treatment for the model
operators $P_B$ and $P_B^t$ in the following. The necessary modifications for
the actual operators $D^{ev}$ and $D^{odd}$ are straightforward and left to
the reader. 

From \cite{Wei} we know that $B$ is essentially selfadjoint if the cone-angles
are less than $2\pi$. We denote by $\bar B$ its closure, which is then a
selfadjoint operator with $\dom \bar B \subset L^2(N)$. Furthermore we know
from \cite{Wei} that under the same assumption on the cone-angles one has
$$
\spec \bar B \cap \textstyle(-\frac{1}{2},\frac{1}{2})=\emptyset,
$$
and hence for a cut-off function $\varphi \in C_0^\infty [0,1)$
$$
\varphi (P_B)_{max} = (P_B)_{min}.
$$
Note that $C_0^{\infty}(N) \subset \dom \bar B$ and $C_0^{\infty}((0,1) \times N) \subset \dom {P_B}_{min}$ are dense with respect to the corresponding graph norms.

Recall from \cite{Wei} that the solution to the homogeneous equation $P_{b}f=0$ is given by
$$
f(r) =r^{-b}f(1),
$$
whereas the solution to the inhomogeneous equation $P_{b}f=g$ is given by
$$
f(r)= r^{-b}\bigl( f(1) + \int_1^r \rho^b g(\rho) d\rho \bigr). 
$$
As in \cite{Wei} we set following \cite{BS}
$$
(T_{b,1}g)(r) = r^{-b}\int_1^r \rho^b g(\rho) d\rho
$$
for $b \in \R$ and
$$
(T_{b,0}g)(r) =  r^{-b}\int_0^r \rho^b g(\rho) d\rho
$$
for $b > -\frac{1}{2}$.

We are mainly interested in the case that $b \not\in(-\frac{1}{2}, \frac{1}{2})$, on which we will focus from now on. We will use $T_{b,1}$ for $b\leq -\frac{1}{2}$ and $T_{b,0}$ for $b \geq \frac{1}{2}$. 

Let now $f \in L^2(0,1)$ and $g \in r^\gamma L^2(0,1)$ for $\gamma \geq 0$, i.e.~$g=r^\gamma \bar{g}$ for $g \in L^2(0,1)$: For $b \leq - \frac{1}{2}$ we get as in Lemma 2.1 in \cite{BS}, resp.~Lemma 4.6 in \cite{Wei}:
\begin{align*}
\vert(T_{b,1}g)(r)\vert \leq  r^{-b} 
& \cdot \left\{ \begin{array}{c@{\quad:\quad}c}
r^{b+\frac{1}{2}+\gamma}\vert 2b+2\gamma+1 \vert^{-\frac{1}{2}} \|\bar g\|_{L^2(0,1)} & b+\gamma <-\frac{1}{2}\\
 \vert \log r \vert^{\frac{1}{2}}\| \bar g\|_{L^2(0,1)} & b+\gamma =-\frac{1}{2}\\
(2b+2\gamma+1)^{-\frac{1}{2}} \|\bar g\|_{L^2(0,1)} & b+\gamma >-\frac{1}{2}
\end{array}\right..
\end{align*}
For $b \geq \frac{1}{2}$ we get
$$
\vert(T_{b,0}g)(r)\vert\leq r^{\frac{1}{2}+\gamma}(2b+2\gamma+1)^{-\frac{1}{2}}\Bigl( \int_0^r |\bar g(\rho)|^2d\rho\Bigr)^{\frac{1}{2}}.
$$
The following is the analogue of Lemma 4.7 in \cite{Wei}:
\begin{lemma}\label{decay_estimates}
Let $f \in L^2(0,1)$ and $g \in r^\gamma L^2(0,1)$ for $\gamma \geq 0$, i.e.~$g=r^\gamma \bar{g}$ for $g \in L^2(0,1)$. Then for $r\in(0,1)$ and with $g=P_{b}f$ we have for $b \geq \frac{1}{2}$
\begin{align*}
|f(r)| &\leq r^{\frac{1}{2}+\gamma}(2b+2\gamma+1)^{-\frac{1}{2}}\Bigl( \int_0^r |\bar g(\rho)|^2d\rho\Bigr)^{\frac{1}{2}}\\
& \leq r^{\frac{1}{2}+\gamma}(2b+2\gamma+1)^{-\frac{1}{2}}\|\bar g\|_{L^2(0,1)},
\end{align*}
and for $b \leq -\frac{1}{2}$
\begin{align*}
|f(r)| \leq   r^{-b}
 \cdot \left\{ \begin{array}{l@{\quad:\quad}l}
|f(1)|+r^{b+\frac{1}{2}+\gamma}\vert 2b+2\gamma+1 \vert^{-\frac{1}{2}} \|\bar g\|_{L^2(0,1)} & b+\gamma <-\frac{1}{2}\\
|f(1)|+ \vert \log r \vert^{\frac{1}{2}}\| \bar g\|_{L^2(0,1)} & b+\gamma =-\frac{1}{2}\\
|f(1)|+(2b+2\gamma+1)^{-\frac{1}{2}} \|\bar g\|_{L^2(0,1)} & b+\gamma >-\frac{1}{2}
\end{array}\right..
\end{align*}
\end{lemma}

\begin{rem}
For $\gamma = 0$ we recover the estimates in Lemma 4.7 in \cite{Wei}.
\end{rem}

As a consequence of these estimates we get the following:

\begin{cor}\label{prekey}
Assume that $B$ is essentially selfadjoint and furthermore
that $\spec \bar B \cap(-\frac{1}{2}, \frac{1}{2}) = \emptyset$. Let $f \in \dom (P_B)_{max}$. Then $f(r) \in L^2(N)$ for all $r \in (0,1)$ and
$$
\|f(r)\|_{L^2(N)} =  
\left\{ \begin{array}{l@{\quad:\quad}l}
O(r^{\frac{1}{2}} |\log r|^{\frac{1}{2}}) & -\frac{1}{2} \in \spec \bar B\\ 
O(r^{\frac{1}{2}}) & -\frac{1}{2}\not\in \spec \bar B
\end{array}\right.
$$
as $r \rightarrow 0$.
\end{cor}

\begin{pf}
Let $\varphi \in C_0^\infty[0,1)$ be a cut-off function satisfying $\varphi\equiv 1$ in a neighbourhood of $0$. Then with $f \in \dom (P_B)_{max}$ we also have $\tilde{f} = \varphi f \in \dom (P_B)_{max}$. According to \cite{Wei} we even have $\tilde f \in \dom (P_B)_{min}$, but we will only make use of $\tilde f(1)=0$ in the following. We set $\tilde g = P_B\tilde f \in L^2((0,1) \times N)$. 
Let $(\psi_b)_{b \in \spec \bar B}$ be an orthonormal system of eigenfunctions of $\bar B$ on $N$. Writing $\tilde f=\sum_b \tilde f_b \otimes \psi_b$ and $\tilde g=\sum_b \tilde g_b \otimes \psi_b$ with $\tilde f_b, \tilde g_b \in L^2(0,1)$ we get $P_b\tilde f_b=\tilde g_b$ for all $b \in \spec \bar B$. The estimates of Lemma \ref{decay_estimates} with $\gamma=0$ (i.e.~Lemma 4.7 in \cite{Wei}) now reduce to
$$
|\tilde{f}_b(r)| \leq
\left\{ \begin{array}{l@{\quad:\quad}l}
r^{\frac{1}{2}}(2b+1)^{-\frac{1}{2}}\|\tilde g_b\|_{L^2(0,1)} & b \geq \frac{1}{2}\\
r^{\frac{1}{2}}|\log r|^{\frac{1}{2}}\|\tilde g_b\|_{L^2(0,1)} & b
 =-\frac{1}{2}\\
r^{\frac{1}{2}}|2b+1|^{-\frac{1}{2}}\|\tilde g_b\|_{L^2(0,1)}
 & b <-\frac{1}{2} \end{array}\right..
$$
By summing over $b \in \spec \bar B$ we get
\begin{align*}
\|\tilde f(r)\|^2_{L^2(N)} &= \sum_b |\tilde f_b(r)|^2\\ 
&\leq 
C \cdot \left\{ \begin{array}{l@{\quad:\quad}l}
r |\log r|  \|\tilde g\|^2_{L^2((0,1) \times N)}& -\frac{1}{2} \in \spec \bar B\\
r \|\tilde g\|^2_{L^2((0,1) \times N)} & -\frac{1}{2} \not\in \spec \bar B
\end{array}\right.,
\end{align*}
hence that $\tilde{f}(r) \in L^2(N)$ and the desired estimate.
\end{pf}\\

The transversal regularity as well as the decay rate can be improved, if $f
\in  \dom (P_B)_{max}(P_B^t)_{max}$:

\begin{lemma}\label{key}
Assume that $B$ is essentially selfadjoint and furthermore
that $\spec \bar B \cap(-\frac{1}{2}, \frac{1}{2}) = \emptyset$. Let $f \in \dom (P_B)_{max}(P_B^t)_{max}$. Then $f(r) \in \dom \bar B$ for all $r \in (0,1)$ and
$$
\|f(r)\|_{\dom \bar B} = O(r^\delta)
$$
as $r \rightarrow 0$, where $\delta \geq 0$ may be any number satisfying both $\delta < \frac{3}{2}$ and $\delta \leq \operatorname{min} \{ b \in \spec \bar B : b \geq \frac{1}{2} \}$.
\end{lemma}

\begin{pf}
Let $\varphi \in C_0^\infty[0,1)$ be a cut-off function as in the proof of
Corollary \ref{prekey}. If $f \in \dom (P_B)_{max}(P_B^t)_{max}$, then we set as above $\tilde f= \varphi f$ and $\tilde g= P^t_B \tilde f$. Again, $\tilde f \in \dom (P^t_B)_{max}$ and $\tilde f(1)=0$. Since $\tilde g = -(\partial_r \varphi) f + \varphi P_B f$, we find that also $\tilde g \in \dom (P_B)_{max}$ and $\tilde{g}(1)=0$. Finally we set $\tilde h=P_B \tilde g \in L^2((0,1) \times N)$. 
Integrating the above estimates, now applied to the equation $P_B\tilde g = \tilde h$, we get that $\tilde g \in r^\gamma L^2((0,1) \times N)$ for any $0 \leq \gamma < 1$. More precisely,
$$
\|r^{-\gamma} \tilde g_b\|^2_{L^2(0,1)} \leq C \cdot
 \left\{ \begin{array}{l@{\quad:\quad}l}
|2b+1|^{-1} \|\tilde h_b\|^2_{L^2(0,1)} & b \neq -\frac{1}{2}\\
\|\tilde h_b\|^2_{L^2(0,1)} & b = -\frac{1}{2}
\end{array}\right..
$$
We will choose $\gamma$ disjoint from the set $-\frac{1}{2}+\spec \bar B$ but
arbitrarily close to $1$ in the following. Then the estimates in Lemma
\ref{decay_estimates} applied to the equation $P^t_B \tilde f = -(P_{-B})
\tilde f = \tilde g$ yield:
\begin{align*}
\|\tilde{f}(r)\|^2_{\dom \bar B} =& \sum_b (1+b^2) |\tilde f_b(r)|^2\\
\leq & r^{2\delta} \sum_b \frac{1+b^2}{|2b+2\gamma+1|} \|r^{-\gamma}\tilde g_b\|^2_{L^2(0,1)}\\
\leq & C r^{2\delta} \Bigl (\sum_{b \neq -\frac{1}{2}} \frac{1+b^2}{|2b+2\gamma+1||2b+1|}\|\tilde h_b\|^2_{L^2(0,1)} + \|\tilde h_{-\frac{1}{2}}\|^2_{L^2(0,1)}\Bigr)\\
\leq & C' r^{2\delta}\|\tilde h\|^2_{L^2((0,1) \times N)}
\end{align*}
for any $\delta \geq 0$ satisfying $\delta < \frac{3}{2}$ and $\delta \leq \operatorname{min} \{ b \in \spec \bar B : b \geq \frac{1}{2}\}$. This proves that $\tilde{f}(r) \in \dom \bar B$ and the corresponding estimate.
\end{pf}\\

The spectrum of $\bar{B}$ has been determined in \cite{Wei}, namely
$$
\spec \bar{B} = \{-1,1\} \cup \Bigl\{ \pm \frac{1}{2} \pm \sqrt{ \frac{1}{4} + \lambda} : \lambda \in \spec \Delta_{N,Fr}, \lambda > 0 \Bigr\}
$$
where $\Delta_{N,Fr}$ is the Friedrichs extension of $\Delta_N=\Delta^0_{d_N}$
on functions. Furthermore $-1$ corresponds to $\ker \Delta_{d_N,Fr}$ on
functions and $1$ corresponds to $\ker \Delta_{d_N,Fr}$ on 2-forms. Note that
the set 
$$
\Bigl\{ \pm \frac{1}{2} \pm \sqrt{ \frac{1}{4} + \lambda} : \lambda \in \spec
  \Delta_{N,Fr}, \lambda > 0 \Bigr\}
$$
does not intersect the interval $[-1,1]$ if
the first positive eigenvalue of $\Delta_{N,Fr}$ is strictly greater than
$2$.  

\begin{cor}\label{cor_key}
Assume that $B$ is essentially selfadjoint and furthermore that $\spec \bar B \cap(-\frac{1}{2}, \frac{1}{2}) = \emptyset$ and that the first positive eigenvalue of $\Delta_{N,Fr}$ is strictly greater than $2$. Let $f \in \dom (P_B)_{max}(P_B^t)_{max}$. 
\begin{enumerate}
\item If $f$ corresponds to an odd form, then
$$
\|f(r)\|_{\dom \bar B} =O(r^\gamma).
$$
\item If $f$ corresponds to a $1$-form, then there exists $\gamma >0$ such that
$$
\|f(r)\|_{\dom \bar B} =O(r^{1+\gamma}).
$$
\end{enumerate}
\end{cor}
\begin{pf}
The first assertion follows directly from Lemma \ref{key}. If $f$
corresponds to a $1$-form, then in the decomposition $f = \sum_b f_b \otimes
\psi_b$ the term corresponding to $b=1$ does not occur (since $f$ does not
involve a $2$-form part on the cross-section). 
\end{pf}\\

A conic differential operator $P$ on $N$ of order $m\geq 0$ acts as a bounded
operator between weighted ``cone'' Sobolev spaces
$$
P : \mathcal{H}^{k+m,\gamma+m}(N) \rightarrow \mathcal{H}^{k,\gamma}(N)\,, k \in
\N_0, \gamma \in \R.
$$
These spaces are nothing but weighted b-Sobolev spaces in the sense of Melrose's
b-calculus, but defined with respect to the measure coming from the cone metric
$g^{TN}$ on $N$; more precisely
\begin{align*}
\mathcal{H}^{k,0}(N)=\{f \in L^2(N) :\,& V_1 \cdot \ldots \cdot V_j f\in
L^2(N) \text{ for all } j \leq k \\
& \text{ and  b-vectorfields } V_i
\}
\end{align*}
and
$$
\mathcal{H}^{k,\gamma}(N)=\rho^\gamma\mathcal{H}^{k,0}(N)
$$
for $\gamma \in \R$, $k \in \N_0$. Here $\rho \in C^\infty(N)$ denotes a
positive function, which close to a cone point equals the distance to that cone point. In particular, we have
$\mathcal{H}^{0,0}(N)=L^2(N)$ and $\mathcal{H}^{0,\gamma}(N)=\rho^\gamma L^2(N)$
for $\gamma \in \R$.

\begin{lemma}
Asssume  that $B$ is essentially selfadjoint and that $0 \not\in\spec \bar B$. Then $\dom \bar B$ is continuously contained in $\mathcal{H}^{1,1}(N)$.
\end{lemma}

\begin{pf}
It follows from \cite{Le} or from \cite{GM} that $\dom \bar B \subset \mathcal{H}^{1,\gamma}(N)$ for any $0 < \gamma <1$. We claim that for such a choice of $\gamma$ the map
$$
B : \mathcal{H}^{1,\gamma}(N) \rightarrow r^{\gamma -1}L^2(N)
$$ 
is invertible. Indeed, since $0 \not\in \spec \bar B$ the $L^2$-kernel of $B$ is trivial, hence so is the kernel of $B$ on $\mathcal{H}^{1,\gamma}(N) \subset L^2(N)$. On the other hand, by duality the cokernel of $B$ on $\mathcal{H}^{1,\gamma}(N)$ is identified with the kernel of $B$ on $\mathcal{H}^{1,1-\gamma}(N)$, which is trivial for the same reason.
Now since $B$ is elliptic, there exists a (generalized) inverse 
$$
G = B^{-1}: r^{\gamma -1}L^2(N) \rightarrow  \mathcal{H}^{1,\gamma}(N),
$$
inducing continuous maps 
$$
G: \mathcal{H}^{k,\mu}(N) \rightarrow \mathcal{H}^{k+1,\mu+1}(N)
$$ 
for all $\mu \in \R$, $k \in \N$. This follows from the existence of a parametrix, cf.~\cite{MeM} or \cite{Sch}, see also \cite{Ma} or \cite{Me}. Now let $f \in \dom \bar B \subset \mathcal{H}^{1,\gamma}(N)$. Then we get $f = GBf \in \mathcal{H}^{1,1}(N)$, since $Bf \in L^2(N)$, and
$$
\|f\|_{\mathcal{H}^{1,1}(N)} \leq C \|Bf\|_{L^2(N)} \leq C \|f\|_{\dom \bar B},
$$
which proves the claim.
\end{pf}

\begin{rem}
It follows from \cite{GM} that actually $\dom \bar B = \mathcal{H}^{1,1}(N)$.
\end{rem}

We summarize what we have achieved so far:

\begin{prop}\label{vertex}
Let $f \in \dom (P_B)_{max}(P_B^t)_{max}$ correspond to a $1$-form. Under the assumptions of Corollary \ref{cor_key} there exists $\gamma > 0$ such that
$$
\|f(r)\|_{\mathcal{H}^{1,1}(N)} = O(r^{1+\gamma})
$$
as $r \rightarrow 0$.
\end{prop}

A direct consequence is the following:

\begin{cor}\label{immediate_cor_vertex}
Let $f \in \dom (P_B)_{max}(P_B^t)_{max}$ correspond to a $1$-form and let $P$ be any first order conic differential operator on $N$. Then under the assumptions of Corollary \ref{cor_key} there exists $\gamma > 0$ such that
$$
\|Pf(r)\|_{L^2(N)} = O(r^{1+\gamma})
$$ 
as $r \rightarrow 0$.
\end{cor}

In the following we mean by a local solution close to a vertex a solution on an open set of the form $(0,\varepsilon) \times N$ for $N$ the link of a vertex $v \in \Sigma$.

\begin{cor}\label{cor_vertex}
Let $\xi$ be a real-valued $1$-form, which is a
local solution of $\Delta_d \xi + 4 \xi=\zeta$ close to a vertex $v$ with $\xi$ in $L^2$, $D\xi$ in $L^2$ and $\zeta$ in $L^2$. If the first positive eigenvalue of $\Delta_{N,Fr}$ is strictly greater than $2$, then there exists $\gamma >0$ such that
$$
\|\xi(r)\|_{L^2(N)}=O(r^\gamma)
$$
and
$$
\| (\nabla_{e_i}\xi)(r)\|_{L^2(N)} = O(r^{\gamma-1})
$$
for $i=2,3$.
\end{cor}

\begin{pf}
If $\xi$ satisfies $\Delta_d \xi + 4 \xi=0$ with $\xi$ in $L^2$, $D\xi$ in $L^2$ and $\zeta$ in $L^2$, then this is equivalent to $D^2 \xi + 4 \xi =\zeta$ with $\xi$ in $L^2$, $D \xi$ in $L^2$ and $\zeta$ in $L^2$. It follows that $D^2\xi$ in $L^2$, i.e.~$\xi \in \dom D_{max}^2$. 

Let $\varphi \in C_0^\infty[0,1)$ be a cut-off function satisfying $\varphi\equiv 1$ in a neighbourhood of $0$. Then by the explicit form of $D^{ev}$ and $D^{odd}$ close to an edge, we get
$$
D( \varphi(r)\xi) = (\partial_r\varphi)(r)\xi + \varphi(r) D \xi
$$ 
and hence that also $\varphi(r)\xi \in \dom D^2_{max}$. Furthermore clearly $\varphi(r)\xi=\xi$ in a neighbourhood of the vertex. Now the result follows from Corollary \ref{immediate_cor_vertex}, resp.~the corresponding statement for $D^{ev}$ in place of $P_B$.
\end{pf}\\

Next we calculate the Hodge-Dirac operator along an edge $e$:
Let $(r,\theta,z)$ be cylindrical coordinates along $e$. Then the hyperbolic metric has the form
$$
g_{hyp} = dr^2 + \sinh(r)^2 d\theta^2 + \cosh(r)^2 dz^2
$$
near $e$. We write
\begin{align*}
\xi^{ev} = \sinh(r)^{-\frac{1}{2}}\cosh(r)^{-\frac{1}{2}} &\bigl ( \varphi_1 +
\varphi_2 \sinh(r) d\theta \wedge dr + \varphi_3 \cosh(r) dz \wedge dr\\
& + \varphi_4 \sinh(r)\cosh(r) d\theta \wedge dz \bigr)
\end{align*}
and
\begin{align*}
\xi^{odd} =\sinh(r)^{-\frac{1}{2}}\cosh(r)^{-\frac{1}{2}}& \bigl (  \varphi_1
dr + \varphi_2 \sinh(r) d\theta + \varphi_3 \cosh(r) dz\\
&+ \varphi_4 \sinh(r)
\cosh(r) dr \wedge d\theta \wedge dz \bigr )
\end{align*}
for even resp.~odd forms. Let $P=D^{ev}$. Then
\begin{align*}
P = \partial_r &- \frac{1}{2} \bigl(\coth(r) + \tanh(r)\bigr)
\\
& + \frac{1}{\sinh(r)} 
\begin{bmatrix}
0 & - \partial_\theta & 0 & 0\\
\partial_\theta & 0 & 0 & 0\\
0 & 0 & \cosh(r) & -\partial_\theta\\
0 & 0 & \partial_\theta & \frac{\cosh(r)^2 + \sinh(r)^2}{\cosh(r)}
\end{bmatrix}\\
&+ \frac{1}{\cosh(r)} 
\begin{bmatrix}
0 & 0 & -\partial_z & 0\\
0 & \sinh(r) & 0 & \partial_z\\
\partial_z & 0 & 0 & 0\\
0 & -\partial_z & 0 & 0
\end{bmatrix}
\end{align*}
acting on quadruples $[\varphi_1,\varphi_2,\varphi_3,\varphi_4]^t$. Hence
\begin{align*}
P_0 &= \partial_r + \frac{1}{r} \begin{bmatrix}
-\frac{1}{2} & -\partial_\theta & 0 &0\\
\partial_\theta & -\frac{1}{2} &0 & 0\\
0&0 & \frac{1}{2}  & -\partial_\theta\\
0&0 & \partial_\theta & \frac{1}{2} 
\end{bmatrix} + 
\begin{bmatrix}
0&0 & -\partial_z & 0\\
0&0 & 0 & \partial_z\\
\partial_z & 0 &0 &0 \\
0 & -\partial_z &0 &0 
\end{bmatrix}\\
&= J \left( \partial_z +\begin{bmatrix}
0&0 & \partial_r & 0\\
0&0 & 0 & -\partial_r\\
-\partial_r & 0 &0 &0 \\
0 & \partial_r &0 &0 
\end{bmatrix}
+ \frac{1}{r}
 \begin{bmatrix}
 0 & 0 & \frac{1}{2} & -\partial_\theta\\
 0 & 0 & -\partial_\theta & -\frac{1}{2}\\
\frac{1}{2}  & \partial_\theta & 0 & 0\\
\partial_\theta & -\frac{1}{2} & 0 & 0
\end{bmatrix}\right)
\end{align*}
with
$$
J=\begin{bmatrix}
0&0 & -1 & 0\\
0&0 & 0 & 1\\
1 & 0 &0 &0 \\
0 & -1 &0 &0 
\end{bmatrix},\quad J^2 = -\id.
$$ 
Similarly for $P^t=D^{odd}$:
\begin{align*}
P^t = -\partial_r  &+ \frac{1}{2} \bigl(\coth(r) + \tanh(r)\bigr)\\
& + \frac{1}{\sinh(r)} 
\begin{bmatrix}
- \frac{\cosh(r)^2 + \sinh(r)^2}{\cosh(r)} & - \partial_\theta & 0 & 0\\
\partial_\theta & -\cosh(r) & 0 & 0\\
0 & 0 & 0 & -\partial_\theta\\
0 & 0 & \partial_\theta & 0
\end{bmatrix}\\
&+ \frac{1}{\cosh(r)} 
\begin{bmatrix}
0 & 0 & -\partial_z & 0\\
0 & 0 & 0 & \partial_z\\
\partial_z & 0 & -\sinh(r) & 0\\
0 & -\partial_z & 0 & 0
\end{bmatrix}
\end{align*}
acting on quadruples $[\varphi_1,\varphi_2,\varphi_3,\varphi_4]^t$. Hence
\begin{align*}
P^t_0  &= -\partial_r + \frac{1}{r}
\begin{bmatrix}
-\frac{1}{2} & -\partial_\theta & 0 &0\\
\partial_\theta & -\frac{1}{2} &0 & 0\\
0&0 & \frac{1}{2}  & -\partial_\theta\\
0&0 & \partial_\theta & \frac{1}{2} 
\end{bmatrix}
 + 
\begin{bmatrix}
0&0 & -\partial_z & 0\\
0&0 & 0 & \partial_z\\
\partial_z & 0 &0 &0 \\
0 & -\partial_z &0 &0 
\end{bmatrix}\\
&=J \left( \partial_z -\begin{bmatrix}
0&0 & \partial_r & 0\\
0&0 & 0 & -\partial_r\\
-\partial_r & 0 &0 &0 \\
0 & \partial_r &0 &0 
\end{bmatrix}
+ \frac{1}{r}
 \begin{bmatrix}
 0 & 0 & \frac{1}{2} & -\partial_\theta\\
 0 & 0 & -\partial_\theta & -\frac{1}{2}\\
\frac{1}{2}  & \partial_\theta & 0 & 0\\
\partial_\theta & -\frac{1}{2} & 0 & 0
\end{bmatrix}\right)
\end{align*}
with $J$ as above.

As before we give a detailed treatment for the model operators $P_0$ and $P_0^t$ and leave the necessary modifications for the actual operators to the reader. 

Let $\varphi \in C_0^\infty[0,1)$ be a cut-off function satisfying $\varphi\equiv 1$ in a neighbourhood of $0$. We consider the operator
$$
A_{\pm}=\pm \begin{bmatrix}
0&0 & \partial_r & 0\\
0&0 & 0 & -\partial_r\\
-\partial_r & 0 &0 &0 \\
0 & \partial_r &0 &0 
\end{bmatrix} + \left( \frac{\varphi(r)}{r} + \frac{1-\varphi(r)}{1-r} \right) \begin{bmatrix}
 0 & 0 & \frac{1}{2} & -\partial_\theta\\
 0 & 0 & -\partial_\theta & -\frac{1}{2}\\
\frac{1}{2}  & \partial_\theta & 0 & 0\\
\partial_\theta & -\frac{1}{2} & 0 & 0
\end{bmatrix}
$$
with domain $C_0^\infty((0,1) \times S^1_\alpha,\R^4) \subset L^2((0,1) \times S^1_\alpha,\R^4)$. Note that $A_{\pm}$ is symmetric and
$$
JA_{\pm}= \pm \partial_r + \left( \frac{\varphi(r)}{r} + \frac{1-\varphi(r)}{1-r} \right) \begin{bmatrix}
-\frac{1}{2} & -\partial_\theta & 0 &0\\
\partial_\theta & -\frac{1}{2} &0 & 0\\
0&0 & \frac{1}{2}  & -\partial_\theta\\
0&0 & \partial_\theta & \frac{1}{2} 
\end{bmatrix}. 
$$
Close to $r=0$ and $r=1$ this operator is of Fuchs type. We calculate its indicial roots at $r=0$, i.e.~we consider
$$
P_B = \partial_r + \frac{1}{r}B
$$ 
with
$$
B=\begin{bmatrix}
-\frac{1}{2} & -\partial_\theta & 0 &0\\
\partial_\theta & -\frac{1}{2} &0 & 0\\
0&0 & \frac{1}{2}  & -\partial_\theta\\
0&0 & \partial_\theta & \frac{1}{2} 
\end{bmatrix}: C^\infty(S^1_\alpha,\R^4) \rightarrow L^2(S^1_\alpha,\R^4). 
$$
An easy calculation shows that
$$
\spec \bar B = \Bigl\{ -\frac{1}{2} \pm \frac{2 \pi n}{\alpha} : n \in \Z \Bigr\} \cup  \Bigl\{ \frac{1}{2} \pm \frac{2 \pi n}{\alpha} : n \in \Z \Bigr\}
$$
and we clearly get the same set of indicial roots at $r=1$. Hence we get by the analysis in \cite{Wei}:

\begin{lemma}
If the cone-angles are less than $2\pi$, then $A_{\pm}$ is essentially selfadjoint. The unique selfadjoint extension $\bar A_\pm$ has discrete spectrum.
\end{lemma}
We may hence decompose
$$
L^2((0,1) \times S^1_\alpha \times (z_0,z_1),\R^4) = \bigoplus_{a \in \spec \bar A_\pm} L^2((z_0,z_1),\R^4)
$$
with $r \in (0,1)$, $\theta \in S^1_\alpha$ and $z \in (z_0,z_1)$. Note that $\dom \bar A_+ = \dom \bar A_-$.

\begin{lemma}\label{transversal_regularity}
Let $f \in \dom (\partial_z + A_\pm)^q_{max}$ for some $q \in \N$. Then
\begin{enumerate}
\item $f(z) \in L^2((0,1) \times S^1_\alpha,\R^4)$ for all $z \in(z_0,z_1)$,
\item $f(z) \in \dom \bar{A}_\pm$ for all $z \in(z_0,z_1)$ if $q \geq 2$,
\item $\partial_zf(z) \in \dom \bar A_\pm$ for all $z \in(z_0,z_1)$ if $q \geq 4$.
\end{enumerate}
\end{lemma}

\begin{pf}
Let $f \in \dom(\partial_z + A_\pm)_{max}$. By multiplying $f$ with a cut-off function $\varphi \in C_0^{\infty}(z_0,z_1)$ we may assume w.l.o.g.~that $f \in \dom (\partial_z +A_\pm)_{min}$ inside $L^2((0,1) \times S^1_\alpha \times (z_0,z_1),\R^4)$. Recall that $C_0^\infty((0,1) \times S^1_\alpha \times (z_0,z_1),\R^4)$ is graph-dense in $\dom (\partial_z +A_\pm)_{min}$. Let now $f \in C_0^\infty((0,1) \times S^1_\alpha \times (z_0,z_1),\R^4)$ and $g = (\partial_z + A_\pm)f$. Let $(\psi_a)_{a \in \spec \bar A_\pm}$ be an orthonormal system of eigenfunctions of $\bar A_\pm$ on $(0,1) \times S^1_\alpha$. Writing $f=\sum_a f_a \otimes \psi_a$ and $g=\sum_a g_a \otimes \psi_a$ we get $\partial_z f_a + a f_a = g_a$ for all $a \in \spec \bar A_\pm$. Solving this ODE with $f_a(z_0) = 0$, resp.~$f_a(z_1) = 0$ we get
\begin{equation}\label{greaterequalzero}
f_a(z) =e^{-az} \int_{z_0}^z e^{w a}g_a(w)dw, 
\end{equation}
resp.
\begin{equation}\label{lesszero}
f_a(z) = - e^{-az} \int_z^{z_1} e^{w a}g_a(w)dw.
\end{equation}
We use (\ref{greaterequalzero}) for $a \in \spec \bar A_\pm$ with $a\geq 0$ and (\ref{lesszero}) for $a \in \spec \bar A_\pm$ with $a< 0$. We thus obtain using the Schwarz inequality
\begin{align}\label{estimates}
\vert f_a(z)\vert \leq 
 \left\{ \begin{array}{c@{\quad,\quad}c}
 (2a)^{-\frac{1}{2}}\|g_a\|_{L^2(z_0,z_1)} & a > 0\\
(z_1-z_0)^{\frac{1}{2}} \|g_a\|_{L^2(z_0,z_1)} & a=0\\
|2a|^{-\frac{1}{2}}\|g_a\|_{L^2(z_0,z_1)} & a < 0
 \end{array}\right.,
\end{align}
i.e.~$|f_a(z)|^2 \leq C(1+a^2)^{-\frac{1}{2}}\|g_a\|^2_{L^2(z_0,z_1)}$.
Summing over $a \in \spec \bar A_\pm$ we get
$$
\|f(z)\|^2_{L^2((0,1) \times S^1_\alpha)} = \sum_a |f_a(z)|^2 \leq C \|g\|^2_{L^2((0,1) \times S^1_\alpha \times (z_0,z_1))}
$$
and by integrating over $z \in (z_0,z_1)$
$$
\|f_a\|^2_{L^2(z_0,z_1)} \leq C (1+a^2)^{-\frac{1}{2}} \|g_a\|^2_{L^2(z_0,z_1)}.
$$
Iterating these estimates we get $f(z) \in \dom \bar A_\pm^{q/2}$ for all $z \in (z_0,z_1)$ if $ f \in \dom(\partial_z+A_{\pm})^q_{max}$, hence assertions 1. and 2.

Finally, if $f \in \dom(\partial_z + A_\pm)_{max}^4$, then $\partial_zf(z) + A_{\pm} f(z) \in \dom \bar A_\pm$ and $f(z) \in \dom \bar A_\pm^2$ for all $z \in (z_0,z_1)$. Therefore $A_\pm f(z) \in \dom \bar A_\pm$ and hence $\partial_zf(z) \in \dom \bar A_\pm$ for all $z \in (z_0,z_1)$. This finishes the proof.
\end{pf}

\begin{lemma}\label{transversal_decay}
Let $f \in \dom \bar A_\pm$. Then $\|f(r)\|_{L^2(S^1_\alpha)}=O(r^{\frac{1}{2}}|\log r|^{\frac{1}{2}})$ as $r \rightarrow 0$.
\end{lemma}

\begin{pf}
This follows from Corollary \ref{prekey}. Note that $-\frac{1}{2} \in \spec \bar B$ in this case. 
\end{pf}\\

In the following we mean by a local solution along an edge a solution on an open set of the form $(0,\varepsilon) \times S^1_\alpha \times (z_0,z_1)$ corresponding to cylindrical coordinates $(r,\theta,z)$ along an edge $e \in \Sigma$.

\begin{cor}\label{cor_edge}
Let $\xi$ be a real-valued $1$-form. If $\xi$ is a local solution of $\Delta_d \xi + 4 \xi=0$ along an edge with $\xi$ in $L^2$ and $D\xi$ in $L^2$, then 
$$
\| \xi(r)\|_{L^2(S^1_\alpha \times (z_0,z_1))} = O(|\log r|^{\frac{1}{2}}) 
$$
and 
$$
\| (\nabla_{e_3}\xi)(r)\|_{L^2(S^1_\alpha \times (z_0,z_1))} = O(|\log r|^{\frac{1}{2}}). 
$$
\end{cor}

\begin{pf}
If $\xi$ satisfies $\Delta_d \xi + 4 \xi=0$ with $\xi$ in $L^2$ and $D\xi$ in $L^2$, then this is equivalent to $D^2 \xi + 4 \xi =0$ with $\xi$ in $L^2$ and $D \xi$ in $L^2$. It follows that $D^2\xi$ in $L^2$, and further by applying powers of $D$ to this equation, that $D^q\xi$ in $L^2$ for any $q \in \N_0$, i.e.~$\xi \in \dom D_{max}^q$ for any $q \in \N_0$. 

Let $\varphi \in C_0^\infty[0,1)$ be a cut-off function satisfying $\varphi\equiv 1$ in a neighbourhood of $0$. Then by the explicit form of $D^{ev}$ and $D^{odd}$ close to an edge, we get
$$
D( \varphi(r)\xi) = (\partial_r\varphi)(r)\xi + \varphi(r) D \xi
$$ 
and hence that also $\varphi(r)\xi \in \dom D^q_{max}$ for any $q \in \N_0$. Furthermore clearly $\varphi(r)\xi=\xi$ in a neighbourhood of the edge. Now the result follows from Lemma \ref{transversal_regularity} and Lemma \ref{transversal_decay}, resp.~the corresponding statements for $P$ in place of $P_0$.  
\end{pf}\\

The estimates given in Corollary \ref{cor_vertex} and in Corollary \ref{cor_edge} will turn out to be sufficient to control the boundary term when integrating by parts later in the main argument.

\section{A vanishing theorem for $L^2$-cohomology}

We begin by reviewing various facts from \cite{HK} and \cite{Wei}, where in case of conflict we prefer to use the notation of \cite{Wei}. 

Let now $(M,g^{TM})$ be the smooth part of a hyperbolic cone-3-manifold and let $\E=\mathfrak{so}(TM) \oplus TM$ denote the bundle of infinitesimal isometries. It carries a canonical flat connection given by
$$
\nabla^\E_X(B,Y) = (\nabla_XB-R(X,Y),\nabla_XY-BX),
$$
where $\nabla$ denotes the Levi-Civita connection on $\E$ and $R$ the Riemannian curvature tensor on $M$. Note that since $M$ has constant sectional curvature $-1$, one has $R(X,Y) = X^\sharp \otimes Y - Y^\sharp\otimes X$, where as usual $X^\sharp$ and $Y^\sharp$ are the 1-forms dual to $X$ and $Y$ with respect to the metric $g^{TM}$.
Let $h^\E$ denote the metric on $\E$ induced by $g^{TM}$. Note however that $\nabla^\E h^\E \neq 0$, whereas of course $\nabla h^\E=0$. Furthermore, the splitting $\E = \mathfrak{so}(TM) \oplus TM$ is $h^\E$-orthogonal and parallel with respect to $\nabla$. 

For $x \in M$ the fibre $\E_x$ may be identified with $\mathfrak{sl}_2(\C)$ and hence has the structure of a Lie algebra. The Lie bracket is given by 
$$
[(A,X),(B,Y)] = ( [A,B]-R(X,Y), AY-BX ),
$$
in particular one has
$$
\nabla^\E_X(B,Y) = \nabla_X(B,Y) + ad(X)(B,Y). 
$$
Note further that $ad(X)$ switches the subbundles $\mathfrak{so}(TM)$ and $TM$: One has $ad(X)(B,0) = (0,-BX)$ and $ad(X)(0,Y)=(-R(X,Y),0)$.

The cross product induces an isomorphism
$$
\times: TM \to \mathfrak{so}(TM), X \mapsto (Y \mapsto X \times Y),
$$
such that we may write $\E = \E' \oplus \E''$ with $\E' \cong \E'' \cong TM$. In fact, since $\mathfrak{sl}_2(\C)$ is a complex Lie algebra and the adjoint representation is $\C$-linear, $\E$ carries a parallel complex structure. Furthermore, $\E'$ is a real form of the complex vector bundle $\E$ and $\E'' = i \E'$. This corresponds to the fact that $\mathfrak{su}(2)$ is a real form of $\mathfrak{sl}_2(\C)$.

Let $\delta^{\E}$ denote the formal adjoint of $d^\E$. Then with respect to a local orthonormal frame $e_1,e_2,e_3$ one has
$$
d^\E = \sum_{i=1}^3 \varepsilon(e^i) \nabla_{e_i}^\E = \sum_{i=1}^3 \varepsilon(e^i) (\nabla_{e_i} + ad(e_i))
$$
and
$$
\delta^{\E} = - \sum_{i=1}^3 \iota(e_i) (\nabla_{e_i}-ad(e_i)).
$$
Let further
$$
D= \sum_{i=1}^3 \varepsilon(e^i) \nabla_{e_i} \quad\text{and}\quad T = \sum_{i=1}^3 \varepsilon(e^i) ad(e_i),  
$$
whose formal adjoints are given by
$$
D^t =  - \sum_{i=1}^3 \iota(e_i) \nabla_{e_i} \quad\text{and}\quad T^t = \sum_{i=1}^3 \iota(e_i) ad(e_i).
$$
Let $\Delta_{d^\E}= d^\E\delta^\E + \delta^\E d^\E$ and let $\Delta_D = DD^t+D^tD$ the associated Laplacians. Let further $H=TT^t+T^tT$. Then one has the following Weitzenb\"ock formula, which is of central importance for our arguments:
\begin{prop}\label{Weitzenboeck}
$\Delta_{d^\E} = \Delta_D + H$.
\end{prop}
The Weitzenb\"ock formula is due to \cite{MM}. The following consequence of the Weitzenb\"ock formula was first observed in \cite{HK}:
\begin{cor}\label{preserve_decomp}
$\Delta_{d^\E}$ preserves the decomposition $\E = \E' \oplus \E''$ in all degrees, i.e.~for $\omega = (\omega',\omega'') \in \Omega^\bullet(M;\E)$ with $\omega' \in \Omega^{\bullet}(M;\E')$ and $\omega'' \in \Omega^{\bullet}(M;\E'')$ one has $\Delta_{d^\E}\omega = (\Delta_{d^\E}\omega',\Delta_{d^\E}\omega'')$.
\end{cor}
\begin{pf}
Since the decomposition $\E = \E' \oplus \E''$ is $h^\E$-orthogonal and parallel with respect to $\nabla$, this is clear for $\Delta_D$. Since $ad(e_i)$ switches the subbundles $\E'$ and $\E''$, the endomorphisms $TT^t$ and $T^tT$ again preserve the above decomposition.
\end{pf}\\
\\
It is easy to see that $TT^t+T^tT$ commutes with the isomorphism induced by the cross product $\times: \Omega^\bullet(M;TM) \rightarrow \Omega^\bullet(M;\mathfrak{so}(TM))$, hence it is enough to compute its action on $\Omega^\bullet(M;TM)$. On 0-forms this was achieved in \cite{HK}:
\begin{lemma}\label{WeitzenboeckremainderI}
$TT^t+T^tT=T^tT=2$ on $\Gamma(M;TM)$.
\end{lemma}
The following is immediate:

\begin{cor}\label{CorWbrI}
$\Delta_{d^\E}(A,X) = (\nabla^t\nabla A, \nabla^t\nabla X)   + 2(A,X)$.
\end{cor}
Note that if $\xi$ denotes the 1-form dual to $X$, then $\nabla^t\nabla\xi + 2 \xi= \Delta_d \xi + 4 \xi$, since by the usual Bochner formula for real-valued 1-forms one has $\Delta_d \xi = \nabla^t\nabla \xi + \Ric \xi =\nabla^t\nabla \xi - 2 \xi$ and for a hyperbolic metric on a 3-manifold $\Ric=-2$.

A 1-form $\eta \in \Omega^1(M;TM)$ may be decomposed into its pure trace, its traceless symmetric and its skew-symmetric part. Clearly $\Delta_D$ and $\nabla^t\nabla$ preserve this decomposition: 
\begin{lemma}\label{WeitzenboeckremainderII}
For $\eta \in \Omega^1(M;TM)$ one has
$$
(DD^t+D^tD) \eta - \nabla^t \nabla \eta = \left\{ 
\begin{array}{c@{\quad:\quad}c}

0 & \eta \text{ pure trace}\\ 
-3\eta & \eta \text{ traceless symmetric}\\
-\eta & \eta \text{ skew-symmetric}
\end{array}
\right..
$$
\end{lemma}

\begin{pf}
If $e_1,e_2,e_3$ is a local orthonormal frame satisfying $\nabla e_i(x)=0$ at a point $x \in M$ and $e^1,e^2,e^3$ the dual coframe, then we compute at $x$
\begin{align*}
DD^t+D^tD &= - \sum_{i,j=1}^3 \bigl(\varepsilon(e^i)\iota(e_j) + \iota(e_i)\varepsilon(e_j)\bigr) \nabla_{e_i}\nabla_{e_j}\\
&= - \sum_{i=1}^3 \nabla_{e_i}\nabla_{e_i} - \sum_{i<j} \bigl(\varepsilon(e^i)\iota(e_j) + \iota(e_i)\varepsilon(e_j)\bigr)R(e_i,e_j)
\end{align*}
Here we have used that $\varepsilon(e^i)\iota(e_j) + \iota(e_j)\varepsilon(e^i)=\delta_{ij}$. Note that we generically use the symbol $R$ for the curvature tensor on any tensor bundle and that
$$
R^{T^*M \otimes TM}(e_i,e_j) = R^{T^*M}(e_i,e_j) \otimes 1 + 1 \otimes R^{TM}(e_i,e_j).
$$
By the usual proof of the Bochner formula on real-valued 1-forms one has
$$
 - \sum_{i<j} \bigl(\varepsilon(e^i)\iota(e_j) + \iota(e_i)\varepsilon(e_j)\bigr)R(e_i,e_j)\otimes 1 = \Ric \otimes 1 =-2,
$$
since $\Ric = -2$ for a hyperbolic metric on a 3-manifold. 
Further, using $R(e_i,e_j)=e^i \otimes e_j - e^j \otimes e_i$ for a hyperbolic metric, we obtain
\begin{align*}
& - \sum_{i<j} \bigl(\varepsilon(e^i)\iota(e_j) + \iota(e_i)\varepsilon(e_j)\bigr)1 \otimes R(e_i,e_j) (e^k \otimes e_l)\\
=&\sum_{i<j} (\delta_{ik}e^j - \delta_{jk}e^i) \otimes (\delta_{il}e_j - \delta_{jl} e_i) = \left\{ 
\begin{array}{c@{\quad:\quad}c}
-e^l \otimes e_k  & k \neq l\\ 
\sum_{i \neq j} e^i \otimes e_i  & k=l=j
\end{array}
\right..
\end{align*}
Evaluating this on the various bits yields the result.
\end{pf}\\
\\
It turns out that also $TT^t+T^tT$ preserves the decomposition of $\Omega^1(M;TM)$ into pure trace, traceless symmetric and skew-symmetric part:

\begin{lemma}\label{WeitzenboeckremainderIII}
For $\eta \in \Omega^1(M;TM)$ one has
$$
(TT^t+T^tT) \eta =  \left\{ 
\begin{array}{c@{\quad:\quad}c}

4 \eta & \eta \text{ pure trace}\\ 
\eta & \eta \text{ traceless symmetric}\\
3 \eta & \eta \text{ skew-symmetric}
\end{array}
\right..
$$
\end{lemma}

\begin{pf}
One easily computes, cf.~also \cite{HK}, that 
$$
TT^t+T^tT = \sum_{i=1}^3 ad(e_i)^2 + \sum_{i,j=1}^3 \varepsilon(e^i)\iota(e_j) ad([e_i,e_j])
$$
with respect to a local orthonormal frame $e_1,e_2,e_3$ and dual coframe $e^1,e^2,e^3$. Now for a vector field $X$
$$
ad(e_i)^2X = R(e_i,X)e_i=(e_i,e_i)X-(X,e_i)e_i,
$$
hence 
$$
\sum_{i=1}^3 ad(e_i)^2X = 2X \quad \text{and} \quad \sum_{i=1}^3 ad(e_i)^2\eta = 2\eta.
$$
Further $[e_i,e_j] = -(e^i \otimes e_j - e^j \otimes e_i)$ and hence
$$
ad([e_i,e_j])X = (X,e_j)e_i-(X,e_i)e_j.
$$
If $\eta= \sum_{i,j=1}^3 \eta_i^j   e^i \otimes e_j$, then
\begin{align*}
\sum_{i,j=1}^3 \varepsilon(e^i)\iota(e_j) ad([e_i,e_j])\eta &= \bigl(\sum_{i=1}^3\eta_i^i\bigr)\bigl(\sum_{i=1}^3 e^i \otimes e_i\bigr) - \sum_{i,j=1}^3 \eta_j^i e^i \otimes e_j\\ &= \tr \eta \cdot \id - \eta^t.
\end{align*}
Evaluating this on the various bits yields the result.
\end{pf}\\
\\
The positivity of $TT^t+T^tT$ on $\Omega^1(M;\E)$, which is originally due to \cite{MM}, follows immediately:

\begin{cor}\label{CorWeitzenboeckremainderIII}
$TT^t+T^tT \geq 1$ on $\Omega^1(M;\E)$.
\end{cor}

As a first step towards the proof of the vanishing theorem we show that both the real and the imaginary part of the $L^2$-harmonic representative of a class in $H^1_{L^2}(M;\E)$ are traceless symmetric, cf.~also \cite{HK}:  

\begin{prop}\label{traceless_symmetric}
Let $\omega \in \Omega^1(M;\E)$ be $L^2$-harmonic. If $\omega=(\omega',\omega'')$ with $\omega' \in \Omega^1(M;\E')$ and $\omega'' \in
\Omega^1(M;\E'')$, then $\omega'$ and $\omega''$ are traceless symmetric.
\end{prop}

\begin{pf}
Note first that $d^\E=D+T$ and $\delta^\E=D^t+T^t$ with $T$ and $T^t$ bundle
endomorphisms which are bounded on $M$ together with all their derivatives. Hence it
is enough to show that a form $\eta \in \Omega^1(M;TM)$ which together with $D\eta, D^t
\eta, DD^t \eta$ and $D^tD\eta$ is in $L^2$ and which satisfies $\Delta_{d^\E}\eta=0$ is in fact
traceless symmetric. 

For $\hat \eta \in \Omega^1(M;TM)$ skew-symmetric let $\xi$ be the real-valued
1-form corresponding to $\hat \eta$. Then  $\Delta_{d^\E}\hat \eta=0$ implies that $\xi$ satisfies the equation
$$
\Delta_d \xi + 4 \xi=0,
$$
cf.~Lemma \ref{WeitzenboeckremainderII}, Lemma \ref{WeitzenboeckremainderIII}
and the remark following Corollary \ref{CorWbrI}. 
We claim that $d\xi, \delta \xi, d\delta \xi$ and $\delta d \xi$ are in $L^2$:
Since $D\hat\eta \in \Omega^2(M;TM)$ is in $L^2$, we find that also
$\star D\hat\eta \in \Omega^1(M;TM)$ is in $L^2$. 
Now a direct calculation shows that  
$
 \delta \xi = - \frac{1}{2}\tr (\star D\hat\eta)
$
and that $d\xi$ is essentially given by the skew-symmetric part of $\star D \hat\eta$. Hence
both of them are in $L^2$. Again by direct calculation we obtain that
$\delta d \xi$ is essentially given by the real-valued 1-form corresponding to
the skew-symmetric part of $DD^t\hat\eta$, hence $\delta d \xi$ and $d \delta \xi$ are
in $L^2$. (Note that $\Delta_d \xi = - 4 \xi$ is in $L^2$.)

Now this implies that $\xi \in \dom d_{max}\delta_{max} + \delta_{max}d_{max}$ and by the
$L^2$-Stokes Theorem $d_{max}\delta_{max} +
\delta_{max}d_{max}$ is a non-negative operator. Note that the $L^2$-Stokes theorem for $\R$-valued forms holds on $M$ since $H^1_{L^2}(N_j;\R)=0$ for all links $N_j$, cf.~\cite{Wei}. We conclude that $\xi = 0$,
hence that $\hat\eta=0$, and in general that the skew-symmetric part of $\eta$ vanishes.

We now analyze the trace part of $\eta$ assuming already that $\eta$ is
symmetric. We find, again using $\Delta_{d^\E}\eta = 0$ together
with Lemma \ref{WeitzenboeckremainderII} and Lemma
\ref{WeitzenboeckremainderIII}, that the trace of $\eta$ satisfies the equation
$$
\Delta_d \tr \eta + 4 \tr \eta = 0.
$$
We claim that $d \tr \eta$ is in $L^2$: Using the symmetry of $\eta$ we find by direct calculation that $d\tr \eta$ is
essentially given by the sum of two terms: The first one is the result of
applying interior multiplication $\iota: \Omega^2(M;TM) \rightarrow
\Omega^1(M)$ to $D\eta$ and the second one is simply the 1-form dual to $D^t
\eta$. 

Hence $\tr \eta \in \dom \delta_{max}d_{max}$ and by the $L^2$-Stokes Theorem $\delta_{max}d_{max}$ is a non-negative operator. We conclude
that $\tr \eta = 0$, hence that $\eta$ is traceless symmetric.
\end{pf}\\

For a vertex $v_j$ let $m_j$ denote the number of edges meeting at $v_j$. Let further $U_\varepsilon(v_j)=B_\varepsilon(v_j) \setminus \Sigma$. Let $N_j$ denote the smooth part of the link of $v_j$.

\begin{lemma}\label{L2cohomlinks}
Let $v_j \in \Sigma$ be a vertex and $N_j$ the smooth part of its link. Then
$$
H^k(N_j;\E) \cong  \left\{ 
\begin{array}{c@{\quad:\quad}c}
0 & k=0,2\\ 
\C^{ 3(m_j-2)} & k=1
\end{array}
\right.
$$
and
$$
H^k_{L^2}(N_j;\E) \cong \left\{ 
\begin{array}{c@{\quad:\quad}c}
0 & k=0,2\\ 
\C^{ 2 (m_j-3)} & k=1
\end{array}
\right.. 
$$
\end{lemma}

\begin{pf}
We drop the index $j$ for convenience. Let $\bar{N} \subset N$ be a compact core with smooth boundary $\partial \bar N$. We use the formula
$$
\chi(\bar N;\E) = \rank(\E) \cdot \chi(\bar N)
$$
(which follows from a Mayer-Vietoris argument applied to a finite good cover of $\bar N$ which trivializes $\E$) and the elementary fact that  $\chi(N) = 2-m$. Since $m \geq 3$ we have $H^0(N;\E)=0$ and, since $N$ is homotopy equivalent to a bouquet of circles, $H^2(N;\E)=0$. Hence it follows that $\dim_\C H^1(N;\E)=3 (m-2)$. To compute the $L^2$-cohomology groups, we use the formula
$$
\chi_{L^2}(N;\E) = \sum_{i=1}^m \chi_{L^2}(U_\varepsilon(p_i);\E) + \chi(\bar{N};\E) - \chi(\partial \bar{N};\E)
$$
from p.~607 in \cite{Ch2}, where $U_\varepsilon(p_i) = B_{\varepsilon}(p_i) \setminus \{p_i\}$ for each cone-point $p_i$. (This formula follows from the Mayer-Vietoris sequence for $L^2$-cohomology, cf.~Lemma 4.3 in \cite {Ch1}, since $\partial \bar N \subset \bar N$ is collared, and the Poincar\'e lemma for collars, cf.~Lemma 3.1 in \cite{Ch1}.) Clearly $H^0_{L^2}(N; \E)=0$ (since already $H^0(N;\E)=0$) and hence, by Poincar\'e duality for $L^2$-cohomology, $H^2_{L^2}(N;\E) = 0$. Note that the $L^2$-Stokes theorem holds for $\E$-valued forms on $N$, since the links (which are just circles) do not have middle-dimensional $\E$-valued $L^2$-cohomology. Furthermore, by Poincar\'e duality for the compact manifold $\partial \bar N$, we have $\chi(\partial \bar N;\E) = 0$. Since $\chi_{L^2}(U_\varepsilon(p_i)) = 1$, we finally obtain $\dim_\C H^1_{L^2}(N;\E) = 2(m-3)$. 
\end{pf}

\begin{lemma}\label{short_exact}
Let $v_j \in \Sigma$ be a vertex and $N_j$ the smooth part of its link. Let furthermore $\bar N_j \subset N_j$ be a compact core with smooth boundary $\partial \bar N_j$. Then there exists a short exact sequence
$$
0 \rightarrow H^1_{L^2}(N_j; \E) \rightarrow H^1(\bar N_j;\E) \rightarrow H^1(\partial \bar N_j;\E) \rightarrow 0.
$$
\end{lemma}

\begin{pf}
We look at a part of the long exact cohomology sequence of the pair
$(\bar{N},\partial \bar{N})$:
$$
\ldots \rightarrow H^1(\bar{N},\partial \bar{N};\E) \overset{q}{\rightarrow}
H^1(\bar{N};\E) \overset{i}{\rightarrow} H^1(\partial \bar{N};\E) \rightarrow \ldots
$$
The map $q$ factors through $H^1_{L^2}(N;\E)$, hence $\dim_\C \im q \leq 2(m-3)$. On the other hand, $\dim_\C H^1(\partial \bar N;\E) =\dim_\C H^0(\partial \bar N;\E)= m$, hence by exactness $\dim_\C \im q \geq 2(m-3)$. Since $\im q \subset \im (H^1_{L^2}(N_j;\E) \rightarrow H^1(N_j;\E))$, this proves the claim.
\end{pf}\\

The preceding lemma is saying that $\E$-valued $L^2$-cohomology may be identified with the space of ordinary cohomology classes in degree $1$, which vanish on the boundary of a compact core. 

We may now compute the $\E$-valued $L^2$-cohomology of singular balls centered at vertices of the singular locus: 

\begin{cor}\label{sing_balls} Let $v_j \in \Sigma$ be a vertex and $U_\varepsilon(v_j)=B_\varepsilon(v_j) \setminus \Sigma$. Then
$$
H^k_{L^2}(U_{\varepsilon}(v_j);\E) \cong \left\{ 
\begin{array}{c@{\quad:\quad}c}
0 & k=0,2,3\\ 
\C^{ 2 (m_j-3)} & k=1
\end{array}
\right.. 
$$
\end{cor}

\begin{pf}
By the Poincar\'e lemma for cones, cf.~Lemma 3.4 in \cite{Ch1}, one has
$$
H^k_{L^2}(U_{\varepsilon}(v_j);\E) \cong \left\{ 
\begin{array}{c@{\quad:\quad}c} 
H^k_{L^2}(N_j;\E) & k=0,1\\
0 & k=2,3
\end{array}
\right.. 
$$
Now the result follows using Lemma \ref{L2cohomlinks}.
\end{pf}\\

Note that due to the presence of middle-dimensional $L^2$-cohomology of the links $N_j$, the $L^2$-Stokes theorem for $\E$-valued forms on $M$ does not hold. This means that $H^\bullet_{max} \neq H^\bullet_{min}$, and in particular that Poincar\'e duality does not hold for $\E$-valued $L^2$-cohomology on $M$.

\begin{lemma}\label{L2cohomlinksII}
Let $e_i \subset \Sigma$ be an edge and $x$ in the interior of $e_i$, let further $N_x$ be the smooth part of its link. Then
$$
H^k_{L^2}(N_x;\E) \cong  \left\{ 
\begin{array}{c@{\quad:\quad}c}
\C  & k=0,2\\ 
0  & k=1
\end{array}
\right.
$$
\end{lemma}

\begin{pf}
We argue in the same way as in the proof of Lemma \ref{L2cohomlinks}. The only difference is that now $m=2$ and $H^0_{L^2}(N_x,\E) \cong H^1_{L^2}(N_x,\E) \cong \C$. 
\end{pf}\\

Finally we compute the $\E$-valued $L^2$-cohomology of singular tubes connecting singular balls centered at the endpoints of an edge:

\begin{cor}\label{L2cohomtubes}
  Let $e$ be an edge with endpoints $v$ and $w$. For $0<\delta\ll\varepsilon$ let $U_\delta(e)=\cup_{x \in e} U_\delta(x) \setminus \Sigma$. Then
$$
H^k_{L^2}(U_\delta(e) \setminus( U_\varepsilon(v) \cup U_\varepsilon(w));\E) \cong \left\{ 
\begin{array}{c@{\quad:\quad}c}
\C & k=0\\ 
0 & k=1,2,3
\end{array}
\right.. 
$$
\end{cor}

\begin{pf}
The result follows from Lemma \ref{L2cohomlinksII}, the Poincar\'e lemma for cones and the Mayer-Vietoris sequence for $L^2$-cohomology, cf.~Lemma 4.3 in \cite{Ch1}. Details are left to the reader.
\end{pf}

\begin{lemma}\label{closedrange}
$H^\bullet_{L^2}(M;\E)$ is finite dimensional and hence the range of $d^\E_{max}$ is closed on $M$.
\end{lemma}

\begin{pf}
This follows from a Mayer-Vietoris argument as in \cite{Ch1} using Corollary \ref{sing_balls} and Corollary \ref{L2cohomtubes}. 
\end{pf}\\

Let $v_j$ be a vertex. Then the holonomy restricted to $N_j$ preserves a point $p_j \in \H^3$, i.e.~$\hol|_{\pi_1(N_j)}$ is conjugate to a representation into $\SU(2)$. Note that the adjoint representation of $\SL_2(\C)$ restricted to $\SU(2)$ acts diagonally on $\sl_2(\C) = \su(2) \oplus i \su(2)$. Therefore we obtain a splitting $\E\vert_{N_j} = \E_j^1 \oplus \E_j^2$ as flat vector bundles, the first summand corresponding to infinitesimal rotations at $p_j$ and the second one corresponding to infinitesimal translations at $p_j$. Furthermore, $Ad \circ \hol | _{\pi_1(N_j)}$ preserves the metric $h^\E_{p_j}$ on 
$$
\E_{p_j} = \mathfrak{so}(T_{p_j}\H^3) \oplus T_{p_j}\H^3,
$$
hence there exists a parallel metric $h_0^\E$ on $\E|_{N_j}$ such that $\E_j^1$ and $\E_j^2$ are orthogonal. We extend the bundles $\E^1_j$ and $\E^2_j$ as well as the metric $h_0^\E$ to $U_\varepsilon(v_j)$ via parallel transport along radial segments.

The metrics $h_0^\E$ and $h^\E$ on $U_\varepsilon(v_j)$ may be compared in the following way, cf.~\cite{Wei}: Let $A$ be the unique field of symmetric endomorphisms of $\E$ on $U_\varepsilon(v_j)$ such that
$$
h^\E(\sigma, \tau) = h_0^\E(A \sigma, \tau)
$$
for $\sigma, \tau \in \Gamma(U_\varepsilon(v_j);\E)$. Let $\delta_0^\E$ denote the formal adjoint of $d^\E$ with respect to the metric $h_0^\E$ (where as usual $\delta^\E$ denotes the formal adjoint of $d^\E$ with respect to $h^\E$). Then according to Lemma 4.2 in \cite{Wei} one has
\begin{equation}\label{divergence}
\delta^\E = \delta_0^\E - A^{-1} \iota(\nabla^\E A),
\end{equation}
where $\iota(\nabla^\E A)$ denotes interior multiplication with the $\End(\E)$-valued 1-form $\nabla^\E A$. Furthermore, according to Lemma 5.7 in \cite{Wei} one has 
\begin{equation}\label{quasi-isometric}
A^{-1}(\nabla^\E A) = - 2 \,ad,
\end{equation} 
which is a bounded $\End(\E)$-valued 1-form.
This implies that $h_0^\E$ and $h^\E$ are quasi-isometric on $U_\varepsilon(v_j)$, cf.~Remark 4.1 in \cite{Wei}. Hence, when computing $H^1_{L^2}(U_\varepsilon(v_j);\E)$, we may replace $h^\E$ by $h_0^\E$. Note however that harmonic forms for these metrics differ since $\delta_0^\E \neq \delta^\E$ in view of (\ref{divergence}) and (\ref{quasi-isometric}).

Let us drop the index $j$ now for convenience. We describe certain standard representatives of $H_{L^2}^1(N;\E)$ and $H^1_{L^2}(U_\varepsilon(v);\E)$ in the following: Let $c^i \in H^1_{L^2}(N;\E^i)$ for $i=1,2$.  
Using Corollary \ref{CorL2Hodge} we may represent $c^i$ by an $L^2$-harmonic form, i.e.~a form $\omega^i \in \Omega^1_{L^2}(N;\E^i)$ satisfying
$$
d_{N}^{\E^i} \omega^i =0 \quad \text{and} \quad \delta_{N}^{\E^i} \omega^i = 0.
$$
Note that the range of $d^{\E^i}$ is closed on $N$ by Lemma \ref{closedrange}. We denote the pull back of this form to $U_{\varepsilon}(v)= (0,\varepsilon) \times N$ by the projection $\pi_N: U_\varepsilon(v) \rightarrow N$ again by $\omega^i$, i.e.~$\omega^i$ is constant in $r$ after identifying the fibers of $\E^i$ using parallel transport along radial segments. Computing $d^\E$ and $\delta_0^\E$ on $U_\varepsilon(v)$ we get 
$$
d^{\E^i} \omega^i =0 \quad \text{and} \quad \delta_0^{\E^i} \omega^i = 0,
$$ 
hence, using (\ref{divergence}) and (\ref{quasi-isometric}), that
\begin{equation}\label{divergenceII}
\delta^\E \omega^i = 2 \iota(ad) \omega^i.
\end{equation}
We have now set up enough notation to state the main result:

\begin{thm}\label{L2vanishing}
Let $M$ be the smooth part of a closed hyperbolic cone-3-manifold $C$ with
cone-angles $\alpha_i \in (0,2\pi)$. Let $\E$ be the flat bundle of
infinitesimal isometries. Let $c \in H^1_{L^2}(M;\E)$ be a class with the
property that for all vertices $v_j$ the following holds:
$$
\left.c \right\vert_{H^1_{L^2}(N_j;\E_j^1)} = 0 \quad\text{or}\quad 
\left.c \right\vert_{H^1_{L^2}(N_j;\E_j^2)} = 0.
$$
Then $c=0$.
\end{thm}
Note that for a trivalent vertex $v_j$ one has $H^1_{L^2}(N_j;\E)=0$ by Lemma \ref{L2cohomlinks}, such that the condition on how $c$ restricts to $N_j$ is empty. In particular, if $\Sigma$ is a trivalent graph, then we get $H^1_{L^2}(M;\E)=0$ by Theorem \ref{L2vanishing}, so we are in situation which is very similar to the one studied in \cite{HK} and \cite{Wei}. On the other hand, in the presence of vertices of valency at least $4$, we will see that indeed $H^1_{L^2}(M;\E) \neq 0$, cf.~Corollary \ref{L2cohom}.

The strategy of the proof of Theorem \ref{L2vanishing} -- which in fact is quite similar to the one in \cite{HK} -- is as follows: For a
class $c \in H^1_{L^2}(M;\E)$ as above let $\omega \in \Omega^1_{L^2}(M;\E)$ denote its $L^2$-harmonic representative given by Corollary \ref{CorL2Hodge} and Lemma \ref{closedrange}, i.e.~$\omega \in \dom d^\E_{max} \cap \dom \delta^\E_{min}$ and $d^\E \omega = \delta^\E \omega = 0$. In order to apply the Bochner method, we wish to justify the integration by parts
$$
0=\int_M |d^\E\omega|^2 + | \delta^\E\omega|^2 = \int_M |D\omega|^2 + |D^t\omega|^2 + (H \omega,\omega).
$$
For this it is sufficient to show that the boundary term
$$
B(r) = - \int_{\partial(M \setminus U_r(\Sigma))} (\id \varotimes h^\E)( \star T\omega \wedge \omega + T^t \omega \wedge \star \omega)
$$
tends to $0$ for a sequence $r \rightarrow 0$.

We do this in 2 steps: First we show that the edges do not give rise to an "ideal" boundary term, i.e.~we remove arbitrarily small balls centered at the vertices and show that integration by parts can be performed on the resulting manifold with boundary (retaining the boundary term of the new boundary components). This is the content of Proposition \ref{edges}. 

Secondly we show that also the vertices do not give rise to an "ideal"
boundary term, i.e.~that integration by parts can be performed on the whole of
$M$. This is the content of Proposition \ref{vertices} and this is where we use the condition on how the class $c$ restricts to the links $N_j$.

\begin{lemma}
Let $\omega=(\omega',\omega'') \in \Omega^1(M;\E)$ be $L^2$-harmonic. Then 
$$
B(r)=2 \int_{\partial(M \setminus U_r(\Sigma))} (\id \varotimes g^{TM}) (\omega' \wedge \omega'').
$$

\end{lemma}

\begin{pf}
From Lemma 2.5 in \cite{HK} it follows that for $\omega=(\omega',\omega'') \in \Omega^1(M;\E)$  with both $\omega'$ and $\omega''$ traceless symmetric one has $T^t \omega = 0$ and $\star T \omega = (\omega'',-\omega')$. Then use Proposition \ref{traceless_symmetric}.
\end{pf}

\begin{prop}\label{edges}
For all $\varepsilon>0$ one has
\begin{align*}
0=\int_{M \setminus \cup_{j=1}^k U_\varepsilon(v_j)} &|D\omega|^2 + |D^t\omega|^2 + (H\omega,\omega)\\
&+ 2 \int_{\partial(M \setminus \cup_{j=1}^k U_\varepsilon(v_j))} (\id \varotimes g^{TM})(\omega' \wedge \omega''),
\end{align*}
i.e.~the edges do not give rise to an "ideal" boundary term.
\end{prop}

\begin{pf}
This is the essence of the argument in \cite{HK}. Let $e$ be an edge with endpoints $v$ and $w$. Let $V_{\varepsilon,r}(e)= \overline{U_r(e)} \setminus (\overline{U_\varepsilon(v)} \cup \overline{U_\varepsilon(w)})$, the closure taken in $M$. Note that $V_{\varepsilon,r}(e)$ is a non-compact manifold with boundary. We show that for a sequence $r \to 0$ the boundary term
$$
B(r) = 2\int_{\partial V_{\varepsilon,r}(e)} (\id \varotimes g^{TM})(\omega' \wedge \omega'') 
$$
tends to zero. We express the boundary term in terms of the orthonormal frame $e_1=\partial/\partial r, e_2=\sinh(r)^{-1}\partial/\partial \theta, e_3= \cosh(r)^{-1} \partial/\partial z$.
We further write $\omega =\sum_{i=1}^3 e^i \otimes \omega_i$, where $e^1,e^2,e^3$ is the coframe dual to $e_1,e_2,e_3$, i.e.~$e^1=dr, e^2= \sinh(r) d \theta, e^3 =\cosh(r) dz$.
Then we obtain
$$
B(r)=  2\int_{\partial V_{\varepsilon,r}(e)} (\omega'(e_2), \omega''(e_3)) - (\omega'(e_3),\omega''(e_2)).
$$
Since by Lemma \ref{L2cohomtubes}, $H^1_{L^2}(U_\delta(e) \setminus( \overline{U_\varepsilon(v)} \cup \overline{U_\varepsilon(w)});\E)=0$ for $0<\delta \ll \varepsilon$, there exists an $L^2$-section $s$ such that $\omega = d^\E s$ on $U_\delta(e) \setminus( \overline{U_\varepsilon(v)} \cup \overline{U_\varepsilon(w)})$. Since $\omega$ is coclosed, we find that $\Delta_{d^\E}s=\delta^\E d^\E s=0$. If we write $s=(s',s'')$ according to the decomposition $\E=\E' \oplus \E''$, then Corollary \ref{preserve_decomp} implies that $\Delta_{d^\E}s'=0$ and $\Delta_{d^\E}s''=0$. If $\xi'$ and $\xi''$ denote the 1-forms corresponding to $s'$ and $s''$, then the remark following Corollary \ref{CorWbrI} implies that $\Delta_d \xi' + 4 \xi'=0$ and $\Delta_d \xi'' + 4 \xi''=0$. Let $\xi$ denote either $\xi'$ or $\xi''$ in the following. Clearly $\xi$ itself is in $L^2$. We claim that also $d\xi$ and $\delta \xi$ are in $L^2$: Since $\omega = d^\E s$ is in $L^2$ and $d^\E=D + T$ with $T$ a bounded $0$-th order operator, we conclude that $\nabla \xi$ is in $L^2$ and hence also $d \xi = \varepsilon \circ \nabla \xi$ and $\delta \xi = - \iota \circ \nabla \xi$. 

We may estimate using the Schwarz inequality
\begin{align*}
\frac{1}{2}|B(r)| &\leq \Bigl(\int_{\partial V_{\varepsilon,r}(e)} |\omega'(e_2)|^2\Bigr)^{1/2} \cdot \Bigl(\int_{\partial V_{\varepsilon,r}(e)} |\omega''(e_3)|^2\Bigr)^{1/2} \\
&+ \Bigl(\int_{\partial V_{\varepsilon,r}(e)} |\omega'(e_3)|^2\Bigr)^{1/2} \cdot \Bigl(\int_{ \partial V_{\varepsilon,r}(e) } |\omega''(e_2)|^2\Bigr)^{1/2}.
\end{align*}
Now $\omega$ in $L^2$ implies that 
$$
\int_{\partial  V_{\varepsilon,r}(e)} |\omega(e_2)|^2 = o(r^{-1} |\log r|^{-1})
$$
for a sequence $r \rightarrow 0$, see Lemma 1.2 in \cite {Ch1}. Since $\xi$ and $D \xi$ are in $L^2$, Corollary \ref{cor_edge} applies. Taking further into account that the volume form on $\partial V_{\varepsilon,r}(e)$ is given by $e^2 \wedge e^3 = \sinh(r)\cosh(r) d\theta \wedge dz$, we get
$$
\int_{\partial  V_{\varepsilon,r}(e)} |\omega(e_3)|^2 = O(r|\log r|). 
$$
For the boundary term we obtain
$$
B(r)=o(1),
$$
again for a sequence $r \rightarrow 0$. This finishes the proof.
\end{pf}\\

We need the following estimate of \cite{MaW} for the first eigenvalue of a spherical cone-surface, which in particular applies to the link of a vertex:
\begin{lemma}\label{ev_estimate}
Let $S$ be a spherical cone-surface with singular locus the collection of points $\{p_1, \ldots, p_m\}$ and cone-angles $\alpha_i \in (0,2\pi)$, $i=1, \ldots,m$. Let $N=S \setminus \{p_1, \ldots, p_m\}$ be the smooth part of $S$ and let $\lambda_1$ be the first (positive) eigenvalue of $\delta_{min}d_{max}$ on functions on $N$. Then $\lambda_1 \geq 2$ with equality $\lambda_1=2$ if and only if $S$ is a spherical suspension.
\end{lemma}
Note that by the $L^2$-Stokes Theorem for $N$ one has $\delta_{min}d_{max}=\delta_{max}d_{min}$, which is the Friedrichs extension of $\Delta_d$ on functions. 

\begin{prop}\label{vertices}
One has
$$
0 = \int_M |D\omega|^2 + |D^t\omega|^2 + (H \omega,\omega),
$$
i.e.~also the vertices do not give rise to an "ideal" boundary term.
\end{prop}

\begin{pf}
Let $v$ be a vertex. We show that for a sequence $r \to 0$ the boundary term
$$
B(r)=2\int_{\partial(M \setminus U_r(v))} (\id \varotimes g^{TM}) (\omega' \wedge \omega'') 
$$
tends to zero. Let $e_1, e_2, e_3$ be a local orthonormal frame with
$e_1=\partial/\partial r$ and let $e^1,e^2,e^3$ be the dual coframe. Note that
$e_2,e_3$ are tangent to $\partial(M \setminus U_r(v))$. If we write $\omega
=\sum_{i=1}^3 e^i \otimes \omega_i$, then we obtain
$$
B(r)=  2\int_{\partial(M \setminus U_r(v))} (\omega'(e_2), \omega''(e_3)) - (\omega'(e_3),\omega''(e_2)).
$$
We now use the assumption on how the class $c$ restricts to $N$ to write
$$
\omega = d^\E s + \left\{ 
\begin{array}{c@{\quad:\quad}c}
\omega^2 & \left.c \right\vert_{H^1_{L^2}(N;\E^1)} =
0\\ 
\omega^1 & \left.c \right\vert_{H^1_{L^2}(N;\E^2)} = 0  
\end{array}
\right. 
$$
on $U_r(v)$ with an $L^2$-section $s$ and a closed form $\omega^i \in
\Omega_{L^2}^1(U_\varepsilon(v);\E^i)$ as described before the statement of Theorem \ref{L2vanishing}.
Since $\omega$ is coclosed, we find that
$
\Delta_{d^\E}s=\delta^\E d^\E s=-\delta^\E \omega^i
$
for either $i=1,2$, and hence, using (\ref{divergenceII}), that
$$
\Delta_{d^\E}s = -2 \iota(ad)\omega^i.
$$
We write $s=(s',s'')$ according to the decomposition $\E=\E' \oplus \E''$ and let $\xi'$ and $\xi''$ be the 1-forms corresponding to $s'$ and $s''$. Then $\Delta_d \xi' + 4 \xi'=\zeta'$ and $\Delta_d \xi'' + 4 \xi''=\zeta''$ with $\xi', \zeta', \xi'',\zeta''$ in $L^2$. Furthermore, as above we have $D\xi', D\xi''$ in $L^2$. 

Let us now assume that $\omega = d^\E s + \omega^1$. The other case is treated in the same way. Then we have that
$$
\lim_{r \rightarrow 0} |(\omega - d^\E s)'-\omega^1|=0 \quad \text{and} \quad
\lim_{r \rightarrow 0} |(\omega - d^\E s)''|=0.
$$
Furthermore $(d^\E s)'=Ds'+Ts''$ and $(d^\E s)'' = Ds''+Ts'$. Therefore
\begin{align*}
\lim_{r \rightarrow 0} B(r) &= 2 \lim_{r \rightarrow 0} \int_{\partial(M \setminus U_r(v))} ((Ds'+Ts''+\omega^1)(e_2),(Ds''+Ts')(e_3))\\
& - 2 \lim_{r \rightarrow 0} \int_{\partial(M \setminus U_r(v))} ((Ds'+Ts''+\omega^1)(e_3),(Ds''+Ts')(e_2)). 
\end{align*}
We estimate the first summand, the second one is treated similarly: We have that $|\omega^1(e_2)|=O(r^{-1})$. Assume that we have a pointwise estimate $s=O(r^\gamma)$ and $\nabla_{e_i}s=O(r^{\gamma-1})$ for $i=2,3$. Since $T$ is a bounded 0-th order operator, we then have $|Ts|=O(r^\gamma)$. Taking into account the volume form $e^2 \wedge e^3 = \sinh(r)^2 dvol_N$ on $\partial(M \setminus U_r(v))$, we get
$$
 \int_{\partial(M \setminus U_r(v))} ((Ds'+Ts'')(e_2),(Ds''+Ts')(e_3))=O(r^{2\gamma})
$$
and
$$
 \int_{\partial(M \setminus U_r(v))} (\omega^1(e_2),(Ds''+Ts')(e_3))=O(r^\gamma).
$$
We conclude that $\lim_{r \rightarrow 0} B(r)=0$ if $\gamma>0$. Now using Corollary \ref{cor_vertex} together with Lemma \ref{ev_estimate} instead of the assumed pointwise estimates we clearly get the same result.
\end{pf}\\

{\slshape Proof of Theorem \ref{L2vanishing}.} By Proposition \ref{vertices} we may integrate by parts. The positivity of the Weitzenb\"ock remainder on $\E$-valued 1-forms, cf.~Corollary \ref{CorWeitzenboeckremainderIII}, yields $\omega=0$.
\hfill$\boxbox$

\section{$L^2$-cohomology and the variety of representations}

Let $\Sigma = e_1 \cup \ldots \cup e_N$, i.e.~$N$ is the number of edges contained in $\Sigma$. Let $k$ be the number of vertices contained in $\Sigma$ and for each vertex $v_j$ let $m_j$ denote the number of edges meeting in $v_j$ (i.e. the valency of the vertex $v_j$). Then one has
$$
2N = \sum_{j=1}^k m_j.
$$
Let $\Sph_j$ denote the link of the $j$-th vertex and $N_j$ its smooth
part. In the following we assume for simplicity that $\Sigma$ is connected and
contains vertices, i.e.~$\Sigma$ is not just a circle.

Let $(r,\theta,z)$ be cylindrical coordinates along an edge $e$ and let $l$ be the length of $e$. We define $\E$-valued forms $\omega_{len}$ and $\omega_{tws}$ on $U_\varepsilon(e)=\cup_{x \in e} U_\varepsilon(x) \setminus \Sigma$ as follows, cf.~\cite{Wei}: Let $\varphi: [0,l] \rightarrow [0,l]$ be a smooth function with
$$
\varphi|_{[0,l/3]} = 0 \text{ and } \varphi|_{[2l/3,l]} = l.
$$
Then we set
$$
\omega_{len} = d \varphi \otimes \sigma_{\partial/\partial z}
$$
and
$$
\omega_{tws} = d \varphi \otimes \sigma_{\partial/\partial \theta}.
$$
We list some properties of $\omega_{len}$ and $\omega_{tws}$ known from \cite{Wei}:
\begin{enumerate}
\item $\omega_{len}$ and $\omega_{tws}$ are closed and in $L^2$.
\item $\omega_{tws} = i \omega_{len}$ with respect to the parallel complex
  structure on $\E$.
\item $\omega_{len}$ and $\omega_{tws}$ infinitesimally do not change the trace
  of the meridian around the edge $e$.
\end{enumerate}
If defined along the $i$-th edge, we denote these forms by $\omega^i_{len}$
and $\omega_{tws}^i$.

\begin{lemma}\label{standardforms} Let $U_{\varepsilon}(\Sigma)= \cup_{x \in
    \Sigma} B_\varepsilon(x) \setminus \Sigma$ and let $N$ denote the number
  of edges and $k$ the number of vertices contained in $\Sigma$. Then 
$$
H^1_{L^2}(U_\varepsilon(\Sigma);\E) = \bigoplus_{i=1}^N \C [\omega^i_{len}]
\oplus \bigoplus_{j=1}^k H^1_{L^2}(U_\varepsilon(v_j);\E)
$$
and in particular
$
\dim_\C H^1_{L^2}(U_\varepsilon(\Sigma);\E) = N + \sum_{j=1}^k 2(m_j - 3)$.
\end{lemma}

\begin{pf}
This follows using the Mayer-Vietoris sequence for $L^2$-cohomology, cf.~Lemma 4.3 in
\cite{Ch1}, together with Corollary \ref{sing_balls} and Corollary \ref{L2cohomtubes}.
\end{pf}

\begin{lemma}\label{injectiveI}
The map $H^1_{L^2}(U_{\varepsilon}(\Sigma);\E) \rightarrow H^1(\partial \bar M;\E)$ is injective.
\end{lemma}

\begin{pf}
This again follows using the Mayer-Vietoris sequence, now for ordinary cohomology, together with Lemma \ref{short_exact} and Lemma \ref{standardforms}.
\end{pf}\\

By the preceding lemma we may identify $H^1_{L^2}(U_{\varepsilon}(\Sigma);\E)$ with a subspace of $H^1(\partial \bar M;\E)$. This subspace is in fact precisely the space of cohomology classes in $H^1(\partial \bar M;\E)$, which vanish on the meridians $\mu_i$ for all $i=1,\ldots,N$, cf.~also Lemma \ref{tangentV}. 

\begin{prop}\label{injectiveII}
The map $H^1_{L^2}(M;\E) \rightarrow H^1(\partial \bar M;\E)$ is
injective.
\end{prop}

\begin{pf}
By Theorem \ref{L2vanishing} a nontrivial class $0 \neq c \in H^1_{L^2}(M;\E)$ restricts
to a nontrivial class in at least one of the groups $H^1_{L^2}(N_j;\E)$, hence
to a nontrivial class already in $H^1_{L^2}(U_{\varepsilon}(\Sigma);\E)$. Now
the result follows from Lemma \ref{injectiveI}.
\end{pf}\\

\begin{cor}\label{injectiveIII}
The map $H^1_{L^2}(M;\E) \rightarrow H^1(\bar M;\E)$ is injective.
\end{cor}

We may hence identify $\E$-valued $L^2$-cohomology on $M$ with a subspace of ordinary cohomology in degree $1$.

\begin{lemma}
$\dim_{\C} H^1(\partial \bar{M};\E) = 6(N-k)$.
\end{lemma}
\begin{pf}
We have
$$
3 \cdot \chi(\partial \bar{M}) = \dim_\C
H^0(\partial \bar{M}; \E) -  \dim_\C
H^1(\partial \bar{M}; \E)
$$
and in our case $H^0(\partial \bar{M};\E)=0$ and $\chi(\partial \bar{M})=2(k-N)$.
\end{pf}
\begin{rem}
\begin{enumerate}
\item In the case that $m_j=3$ for all $j$, one has $2N=3k$ and hence $6(N-k)=2N$.
\item In the general case one has $6(N-k)= 2N+\sum_{j=1}^k 2(m_j-3)$. 
\item If $\mathcal{T}_{0,m}$ denotes the Teichm\"uller space of the $m$-times punctured sphere, then one has $\dim_{\C}\mathcal{T}_{0,m}=m-3$.
\end{enumerate}
\end{rem}
\begin{prop}\label{half_dim}
The natural map $i : H^1(\bar{M};\E) \rightarrow H^1(\partial \bar{M};\E)$ is injective and $\dim_{\C} H^1(\bar{M};\E) = \frac{1}{2}\dim_{\C} H^1(\partial \bar{M};\E
) = 3(N-k)$.
\end{prop}
\begin{pf}
We look at a part of the long exact cohomology sequence of the pair
$(\bar{M},\partial \bar{M})$:
$$
\ldots \rightarrow H^1(\bar{M},\partial \bar{M};\E) \overset{q}{\rightarrow}
H^1(\bar{M};\E) \overset{i}{\rightarrow} H^1(\partial \bar{M};\E) \rightarrow \ldots
$$
Let now $c \in H^1(\bar M;\E)$ with $i(c)=0$. By exactness at
$H^1(\bar{M};\E)$ there exists $b \in
H^1(\bar{M},\partial \bar{M};\E)$ with $q(b)=c$. Since $q$ factors through
$H^1_{L^2}(M;\E)$, using Proposition \ref{injectiveII} we conclude $c=0$. This
proves injectivity of $i$. 

Furthermore, Poincar\'e duality yields the short exact sequence
$$
0 \rightarrow \im q \rightarrow H^1(\bar{M};\E)
\rightarrow H^1(\partial \bar{M};\E) \rightarrow H^1(\bar{M},\E)^*
\rightarrow \im q^*
\rightarrow 0,
$$ 
see \cite{HK} or \cite{Wei} for details. The result follows taking into account that $q=0$ since $i$ is injective.
\end{pf}
\begin{rem}
\begin{enumerate}
\item In the case that $m_j=3$ for all $j$, one has $3(N-k)=N$.
\item In the general case one has $3(N-k)= N+\sum_{j=1}^k (m_j-3)$.
\end{enumerate}
\end{rem}

Let now $v \in \Sigma$ be a vertex and $N$ the smooth part of its link $\Sph$. Let further $p_1, \ldots ,p_m \in \Sph$ be the cone points. Hence $N = \Sph \setminus \{p_1, \ldots, p_m\}$
is homeomorphic to the $m$-times punctured sphere $S^2 \setminus \{p_1, \ldots, p_m\}$ and
$$
\pi_1 N = \langle \gamma_1, \ldots, \gamma_m \,\vert\, \gamma_1 \cdot \ldots \cdot \gamma_m \rangle,
$$
the free group of rank $m-1$. Here the $\gamma_i$ are the obvious loops around the punctures $p_i$. It follows that
$$
R(\pi_1 N, \SL_2(\C))=\{ (A_1, \ldots, A_m) \in \SL_2(\C)^m:A_1\cdot \ldots \cdot A_m=1\}.
$$
Clearly the map $f: \SL_2(\C)^m \rightarrow \SL_2(\C), (A_1, \ldots,A_m) \mapsto A_1 \cdot \ldots \cdot A_m$ is a submersion, such that $R(\pi_1 N, \SL_2(\C))=f^{-1}(1) \subset \SL_2(\C)^m$ is a smooth submanifold of complex dimension $3(m-1)$. Furthermore, as done in \cite{Wei} for the case $m=3$, one shows that the map
$$
t_\gamma: R(\pi_1 N, \SL_2(\C)) \rightarrow \C^m, \rho \mapsto (t_{\gamma_1}(\rho), \ldots, t_{\gamma_m}(\rho))
$$
with $t_{\gamma_i}(\rho):=\tr \rho(\gamma_i)$ is a submersion at the holonomy representation, which we will also denote by $\rho_0$ in the following. Equivalently this means that the differentials $dt_{\gamma_1}, \ldots, dt_{\gamma_m}$ are $\C$-linearly independent in $T^*_{\rho_0}R(\pi_1N,\SL_2(\C))$.

We consider the map
$$
t_\mu: R(\pi_1 \partial \bar M, \SL_2(\C)) \rightarrow \C^N, \rho \mapsto (t_{\mu_1}(\rho), \ldots, t_{\mu_N}(\rho)),
$$
where the $\mu_i$ are the meridian loops around the edges $e_i$.

\begin{lemma}\label{Rsmooth}
The representation $\rho_0$ is a smooth point in $R(\pi_1\partial \bar{M}, \SL_2(\C))$, furthermore the differentials $\{dt_{\mu_1}, \ldots ,dt_{\mu_N}\}$ are $\C$-linearly independent in $T^*_{\rho_0} R(\pi_1\partial \bar{M}, \SL_2(\C))$. The local $\C$-dimension of $R(\pi_1\partial \bar{M}, \SL_2(\C))$ at the representation $\rho_0$ equals $2N + \sum_{j=1}^k 2(m_j-3) + 3$.
\end{lemma}

\begin{pf}
Using the above facts about $R(\pi_1N,\SL_2(\C))$ the glueing procedure goes through as described in \cite {Wei} for the case $m_j=3$ for $j=1,\ldots,k$. A careful dimension count yields the formula for $\dim_\C R(\pi_1\partial \bar{M}, \SL_2(\C))$ at $\rho_0$. Details are left to the reader.
\end{pf}\\

Since the holonomy representation of a hyperbolic cone-manifold structure is irreducible (see \cite{HK} or \cite{Wei} in the presence of vertices) and the action of $\SL_2(\C)$ on the irreducible part of $R(\pi_1 \partial \bar M, \SL_2(\C))$ is proper (see for example Lemma 6.24 in \cite{Wei}), we obtain as in Corollary 6.25 in \cite{Wei} the following statement:

\begin{lemma}\label{Xsmooth}
The equivalence class $\chi_0$ of the representation $\rho_0$ is a smooth point in $X(\pi_1\partial \bar{M}, \SL_2(\C))$. The local $\C$-dimension of $X(\pi_1\partial \bar{M}, \SL_2(\C))$ at $\chi_0$ equals $2N + \sum_{j=1}^k 2(m_j-3)$. The tangent space $T_{\chi_0}X(\pi_1\partial \bar{M}, \SL_2(\C))$ may be identified with $H^1(\partial \bar M; \E)$.
\end{lemma}

Since the traces are constant on the orbits of the action of $\SL_2(\C)$ on $R(\pi_1 \partial \bar M,\SL_2(\C))$, the differentials
$\{dt_{\mu_1}, \ldots ,dt_{\mu_N}\}$ remain $\C$-linearly independent in
$T^*_{\chi_0} X(\pi_1\partial \bar{M}, \SL_2(\C))=H^1(\partial \bar M;\E)^*$. Hence
the level set
$$
V =\{ t_{\mu_1}=
t_{\mu_1}(\chi_0), \ldots, t_{\mu_N}=t_{\mu_N}(\chi_0 )\}
$$
is locally at $\chi_0$ a smooth submanifold in $X(\pi_1\partial \bar{M}, \SL_2(\C))$ of $\C$-dimension $N + \sum_{j=1}^k 2(m_j-3)$ and with
tangent space
$$
T_{\chi_0}V= \{ dt_{\mu_1}= \ldots = dt_{\mu_N}=0\}.
$$
\begin{lemma}\label{tangentV}
$T_{\chi_0}V=H^1_{L^2}(U_\varepsilon(\Sigma);\E)$.
\end{lemma}

\begin{pf}
By Lemma \ref{short_exact} and Lemma \ref{standardforms} and the fact that $dt_{\mu_{i'}}([\omega^i_{len}])=0$ we get that $H^1_{L^2}(U_\varepsilon(\Sigma);\E) \subset T_{\chi_0}V$. Now a dimension argument yields the result.
\end{pf}\\

As in \cite{Wei}, using a construction of M.~Kapovich, cf.~Lemma 8.46 in \cite{Kap}, and the irreducibility of $\rho_0$, we get:

\begin{lemma}
The equivalence class $\chi_0$ of the representation $\rho_0$ is a smooth point in $X(\pi_1 \bar M,\SL_2(\C))$. Its tangent space $T_{\chi_0} X(\pi_1 \bar M,\SL_2(\C))$ may be identified with $H^1(\bar M;\E)$.
\end{lemma}
Together with Proposition \ref{half_dim} this yields the following statement:
\begin{cor}\label{deformationspace}
The local $\C$-dimension of $X(\pi_1 \bar{M}, \SL_2(\C))$ at $\chi_0$ equals $N + \sum_{j=1}^k (m_j-3)$.
\end{cor}
Using the holonomy theorem of Ehresmann-Thurston we obtain:
\begin{cor}
The deformation space $\operatorname{Def}(M)$ of hyperbolic structures on $M$ is locally homeomorphic to $\C^{N + \sum_{j=1}^k (m_j-3)}$.
\end{cor}

Recall at this point that $\operatorname{Def}(M)$ is the space of all deformations into incomplete hyperbolic structures, which are not necessarily cone-manifold structures. The space of cone-manifold structures $C_{-1}(X,\Sigma) \subset \operatorname{Def}(M)$ of fixed topological type $(X,\Sigma)$ is a proper subspace.

\begin{prop}\label{traces_indep}
The differentials $\{dt_{\mu_1}, \ldots ,dt_{\mu_N}\}$ are $\C$-linearly
independent already in $T^*_{\chi_0} X(\pi_1 \bar{M}, \SL_2(\C))=H^1(\bar M;\E)^*$.
\end{prop}

\begin{pf}
By what has been said above, it is enough to show that the subspaces
$H^1_{L^2}(U_\varepsilon(\Sigma);\E), H^1(\bar M;\E) \subset H^1(\partial \bar
M;\E)$ meet transversally in the sense that
$$
H^1_{L^2}(U_\varepsilon(\Sigma);\E) + H^1(\bar M;\E) =H^1(\partial \bar M;\E).
$$
Recall that $$\dim_\C H^1_{L^2}(U_\varepsilon(\Sigma);\E) = N + \sum_{j=1}^k
2(m_j-3)$$ and $$\dim_\C H^1(\bar M;\E)= N + \sum_{j=1}^k (m_j-3),$$ whereas $\dim_\C
H^1(\partial \bar M;\E) = 2N + \sum_{j=1}^k 2(m_j-3)$. It follows that
\begin{equation}\label{lowerbound}
\dim_\C H^1_{L^2}(U_\varepsilon(\Sigma);\E) \cap H^1(\bar M;\E) \geq \sum_{j=1}^k (m_j-3).
\end{equation} 
On the other hand, by Theorem \ref{L2vanishing} the space $
H^1_{L^2}(U_\varepsilon(\Sigma);\E) \cap H^1(\bar M;\E)$ is contained in
$\bigoplus_{j=1}^k H^1_{L^2}(U_\varepsilon(v_j);\E)$ and intersects
$\bigoplus_{j=1}^k H^1_{L^2}(U_\varepsilon(v_j);\E_j^1)$ trivially. It follows that
\begin{equation}\label{upperbound}
\dim_\C H^1_{L^2}(U_\varepsilon(\Sigma);\E) \cap H^1(\bar M;\E) \leq \sum_{j=1}^k (m_j-3).
\end{equation} 
Now we compute
$
\dim_\C  (H^1_{L^2}(U_\varepsilon(\Sigma);\E) + H^1(\bar M;\E)) = 2N +
\sum_{j=1}^k 2(m_j-3) = \dim_\C H^1(\partial \bar M;\E)
$
which proves the assertion.
\end{pf}

\begin{cor}\label{L2cohom}
$\dim_\C H^1_{L^2}(M;\E) = \sum_{j=1}^k (m_j-3)$.
\end{cor}

\begin{pf}
Clearly $H^1_{L^2}(U_\varepsilon(\Sigma);\E) \cap H^1(\bar M;\E) = \im(
H^1_{L^2}(M;\E) \to H^1(\bar M;\E))$ and by Corollary \ref{injectiveIII}
this map is injective. Hence we get 
$$
H^1_{L^2}(U_\varepsilon(\Sigma);\E) \cap H^1(\bar M;\E)=H^1_{L^2}(M;\E).
$$
From \eqref{lowerbound} and \eqref{upperbound} we obtain the result.  
\end{pf}\\  

Consider the map
$$
t_{\mu}=(t_{\mu_1}, \ldots, t_{\mu_N}): X(\pi_1 \bar{M}, \SL_2(\C)) \rightarrow \C^N
$$
and its differential $(dt_\mu)_{\chi_0}: H^1(M;\E) \rightarrow \C^N$. Proposition \ref{traces_indep} says that $t_\mu$ is a submersion at $\chi_0$, i.e.~that $t_{\mu_1}, \ldots, t_{\mu_N}$ is part of a local coordinate system. Note that $\dim_\C X(\pi_1 \bar M,\SL_2(\C))-N = \sum_{j=1}^k (m_j-3)$, which is the number of missing coordinates. It remains to construct these missing coordinates. We will return to this question in \cite {MoW}.

\begin{cor}\label{kerneltraces}
$\ker (dt_\mu)_{\chi_0} = H^1_{L^2}(M;\E)$.
\end{cor}

\begin{pf}
This follows from Lemma \ref{tangentV}.
\end{pf}\\

Let $X_0(\pi_1 \partial \bar M, \SL_2(\C))$ denote the space of equivalence classes of representations $\rho: \pi_1 \partial \bar M \rightarrow \SL_2(\C)$ such that for all vertices $v_j$ the restriction $\rho\vert_{N_j}$ fixes a point $p_j \in \H^3$, i.e.~$\rho\vert_{N_j}$ is conjugate to a representation into $\SU(2)$. 

\begin{lemma}
$\chi_0$ is a smooth point in $X_0(\pi_1 \partial \bar M, \SL_2(\C))$ and furthermore
$
\dim_\R X_0(\pi_1 \partial \bar M, \SL_2(\C))= 3N + \sum_{j=1}^k 2(m_j - 3)
$
at $\chi_0$.
\end{lemma}

\begin{pf}
This is done using the same constructions as in Lemma \ref{Rsmooth} and in Lemma \ref{Xsmooth}. Details are left to the reader.
\end{pf}\\

The preceding lemma relates nicely to the deformation space of a cone-tube: The $3N$ real parameters correspond to $3$ real parameters for each edge, namely the cone-angle, the length and the twist of the edge. 
Furthermore, the $2(m_j-3)=\dim_{\R}\mathcal{T}_{0,m_j}$ real parameters for each vertex correspond to the conformal part of the deformation space of the link, cf.~\cite{Tro} and \cite{LT}, see also \cite{MaW}.

Continuing our main argument, we observe that 
\begin{align*}
\dim_\R X_0(\pi_1 \partial \bar M, \SL_2(\C)) &+ \dim_\R X(\pi_1 \bar M, \SL_2(\C)) = 5N + \sum_{j=1}^k 4(m_j-3)\\ &= \dim_\R X(\pi_1 \bar M, \SL_2(\C)) + N
\end{align*}
such that $X_0(\pi_1 \partial \bar M, \SL_2(\C))$ and $X(\pi_1 \bar M, \SL_2(\C))$ meet transversally at $\chi_0$ and the intersection
$$
X_0(\pi_1 \bar M, \SL_2(\C)):=X_0(\pi_1 \partial \bar M, \SL_2(\C)) \cap X(\pi_1 \bar M, \SL_2(\C)).
$$
is locally a smooth submanifold with $\dim_\R X_0(\pi_1 \bar M, \SL_2(\C)) \geq N$ at $\chi_0$.

\begin{thm}\label{premain}
$\dim_\R X_0(\pi_1 \bar M, \SL_2(\C))= N$ at $\chi_0$ and the map
$$
t_{\mu}=(t_{\mu_1}, \ldots, t_{\mu_N}): X_0(\pi_1 \bar M, \SL_2(\C)) \rightarrow \R^N 
$$
is a local diffeomorphism at $\chi_0$.
\end{thm}

\begin{pf}
We claim that $(dt_\mu)_{\chi_0}$ is injective. Indeed, let $c \in H^1(M;\E)$ be a class with $dt_\mu(c)=0$. Then we get using Corollary \ref{kerneltraces} that $c \in H^1_{L^2}(M;\E)$ considered as a subspace in $H^1(M;\E)$. Now since $c$ is also tangent to $X_0(\pi_1 \partial \bar M, \SL_2(\C))$ we obtain that $c|_{H^1_{L^2}(N_j;\E_j^2)}=0$ for all vertices $v_j$. Hence Theorem \ref{L2vanishing} applies to yield $c=0$. 

Finally, injectivity of $(dt_\mu)_{\chi_0}$ yields $\dim_\R X_0(\pi_1 \bar M, \SL_2(\C)) \leq N$, hence $\dim_\R X_0(\pi_1 \bar M, \SL_2(\C)) = N$ and further that $t_\mu$ is a local diffeomorphism at $\chi_0$.
\end{pf}\\

As a consequence we obtain our main result:

\begin{thm}[Local Rigidity]\label{main}
Let $X$ be a hyperbolic cone-3-manifold with cone-angles less than $2\pi$. Then the map 
$$
\alpha=(\alpha_1, \ldots, \alpha_N): C_{-1}(X,\Sigma) \rightarrow \R_+^N
$$ 
is a local homeomorphism at the given structure.
\end{thm}

\begin{pf}
We just have to apply Theorem \ref{premain} and the Ehresmann-Thurston holonomy theorem together with the usual relation between the trace of the meridians and the cone-angles. 
\end{pf}


\end{document}